\documentclass[11pt,a4paper]{amsart}
\usepackage{amssymb}
\usepackage[T1]{fontenc}
\usepackage{Alegreya}
\usepackage[libertine]{newtxmath}

\usepackage{mathrsfs}
\usepackage{tikz-cd}
\usepackage[headings]{fullpage}
\usepackage{paralist}
\usepackage{xspace}

\usepackage{color}
\usepackage[normalem]{ulem}

\usepackage[colorlinks=true,linkcolor=blue, citecolor=blue, urlcolor=blue,%
pagebackref=true]{hyperref}

\theoremstyle{plain} \newtheorem{theorem}[equation]{Theorem}
\newtheorem{lemma}[equation]{Lemma}
\newtheorem{corollary}[equation]{Corollary}
\newtheorem{proposition}[equation]{Proposition}

\theoremstyle{definition}

\newtheorem{remark}[equation]{Remark} \newenvironment{remarkbox}[1][]{%
\begin{remark}[#1] \pushQED{\qed}}{\popQED \end{remark}}

\newtheorem{example}[equation]{Example} 

\newtheorem{definition}[equation]{Definition}

\newtheorem{notation}[equation]{Notation}

\newtheorem{discussion}[equation]{Discussion}
\newenvironment{discussionbox}[1][]{%
\begin{discussion}[#1]\pushQED{\qed}}{\popQED \end{discussion}}

\newtheorem{observation}[equation]{Observation}

\newtheorem{construction}[equation]{Construction}

\newtheorem{setup}[equation]{Setup}
 \newcounter{step}
\newcommand{\fraka}{{\mathfrak a}}
\newcommand{\frakb}{{\mathfrak b}}
\newcommand{\calI}{\mathcal I}

\newcommand{\frakM}{{\mathfrak M}}
\newcommand{\frakm}{{\mathfrak m}}

\newcommand{\frakn}{{\mathfrak n}}
\newcommand{\frakp}{{\mathfrak p}}
\newcommand{\bfx}{\mathbf x}

\newcommand{\naturals}{\mathbb{N}}
\newcommand{\ints}{\mathbb{Z}}
\def\to{\longrightarrow}
\newcommand{\Rees}{\mathscr R}
\DeclareMathOperator{\Hom}{Hom}
\DeclareMathOperator{\Spec}{Spec}
\DeclareMathOperator{\Proj}{Proj}
\DeclareMathOperator{\Min}{Min}
\DeclareMathOperator{\depth}{depth}
\DeclareMathOperator{\homology}{H}
\newcommand{\define}[1]{\emph{#1}}
\newcommand{\minus}{\ensuremath{\smallsetminus}}
\DeclareMathOperator{\image}{Im}
\DeclareMathOperator{\Cech}{\check C}
\newcommand{\tCech}{\v Cech\,}
\DeclareMathOperator{\socle}{soc}
\DeclareMathOperator{\ord}{ord}

\newif\ifreadkumminibib
\readkumminibibfalse

\begin{document}
\title{Blow-up rings and $F$-rationality}
\author{Nirmal Kotal}
\address{Chennai Mathematical Institute, Siruseri, Tamilnadu 603103. India}
\email{nirmal@cmi.ac.in}
\author{Manoj Kummini}
\address{Chennai Mathematical Institute, Siruseri, Tamilnadu 603103. India}
\email{mkummini@cmi.ac.in}
\thanks{Manoj Kummini was partially supported by the grant MTR/2017/00006
from Science and Engineering Research Board, India. Both authors
were partly supported by an Infosys Foundation fellowship.}
\subjclass{13A35; 13A30}
\keywords{$F$-rationality, test ideals, Rees algebras}
\begin{abstract}
In this paper, we prove some sufficient conditions for Cohen-Macaulay
normal Rees algebras to be $F$-rational.
Let $(R,\mathfrak{m})$ be a Gorenstein normal local domain of dimension
$d\geq 2$ and of characteristic $p > 0$.
Let $I$ be an ideal generated by a system of parameters.
Our first set of results give conditions on the test ideals $\tau(I^n)$, $n
\geq 1$ which would imply that the normalization of the Rees algebra
$R[It]$ is $F$-rational.
Another sufficient condition is that the socle of
$\homology_{\overline{G}_+}^d(\overline{G})$ (where $\overline{G}$ is the
associated graded ring for the integral closure filtration) is entirely in
degree $-1$, if $R$ is $F$-rational (but not necessarily Gorenstein).
Then we show that if $R$ is a hypersurface of degree $2$ or is
three-dimensional and $F$-rational and $\Proj (R[\frakm t ])$ is
$F$-rational, then $R[\frakm t ]$ is $F$-rational.
\end{abstract}

\maketitle

\section{Introduction}
\numberwithin{equation}{section}

Let $(R, \frakm)$ be a $d$-dimensional excellent local
domain of prime characteristic $p>0$, where $d \geq 2$.
Let $I$ be an $\frakm$-primary ideal.
Write $\Rees(I)$ for the Rees algebra $\oplus_{n \in \naturals}I^n$ and
$\overline{\Rees(I)}$ for its normalization.
In this paper, we prove some sufficient conditions for
$\overline{\Rees(I )}$ to be $F$-rational. These are motivated
by results of
N.~Hara, K.-i.~Watanabe and
K.-i.~Yoshida~\cite{HaraWatanabeYoshidaFrationality2002},
Hara and Yoshida~\cite{HaraYoshidaGenTightClMultIdeals2003},
M.~Koley and the second author~\cite{KoleyKumminiFrationality2017}
and the analogous results of
E.~Hyry~\cite{HyryBlowupRingsRationalSings1999}
in characteristic zero.

Our first result is the following converse
to~\cite[Theorem~5.1]{HaraYoshidaGenTightClMultIdeals2003}.
For an $R$-ideal $\fraka$, $\tau(\fraka)$ is the \define{test ideal} of
$\fraka$~\cite[1.9]{HaraYoshidaGenTightClMultIdeals2003}.
\begin{theorem}
\label{theorem:frationalityifftestideal}
Let $(R,\mathfrak{m})$ be a Gorenstein normal local domain of dimension
$d\geq 2$ and of characteristic $p > 0$.
Let $I$ be an $R$-ideal generated by a system of  parameters such that
$\overline{\Rees(I)}$ is Cohen-Macaulay.
If $\tau(I^n)=I^n:\overline{I^{d-1}}$ for all integers $1 \leq n \leq d-1$,
then $\overline{ \Rees(I)}$ is $F$-rational.
\end{theorem}

Using Lemma~\ref{lemma:vanishingzerostar}, we immediately obtain the
following:
\begin{corollary}
Write $\overline{\mathfrak{M}}$ for the unique
homogeneous maximal ideal of $\overline{\Rees(I )}$.
If
\[
\left[ 0^{*}_{\homology^{d+1}_{\overline{\mathfrak{M}}} (\overline{\Rees(I)})}
\right]_{-n} = 0
\]
for all $1 \leq n \leq d-1$, then
$\overline{\Rees(I )}$ is $F$-rational.
\end{corollary}

As a corollary of some of the arguments that go into the proof of
Theorem~\ref{theorem:frationalityifftestideal}, we get the following
proposition.

\begin{proposition}
\label{proposition:gorG}
Let $(R,\mathfrak{m})$ be Gorenstein normal local domain of dimension
$d\geq 2$ and of characteristic $p> 0$.
Let $I$ be an $R$-ideal generated by a system of  parameters such that
$\overline{\Rees(I)}$ is Cohen-Macaulay and
the associated graded ring
$\overline{G} := \oplus_{n \geq 0}
\frac{\overline{I^n}}{\overline{I^{n+1}}}$
is Gorenstein.
Write $a$ for the $a$-invariant of $\overline{G}$.

\begin{enumerate}

\item
\label{proposition:gorG:fwd}
If $\tau(I^{-a-1}) = R$, then $\overline{\Rees(I)}$ is $F$-rational.

\item
\label{proposition:gorG:conv}
If $\overline{\Rees(I)}$ is $F$-rational and $a \leq -2$, then
$\tau(I^n) = R$ for all $0 \leq n \leq -a-1$.
\end{enumerate}
\end{proposition}

The next proposition overlaps with
Proposition~\ref{proposition:gorG}~\eqref{proposition:gorG:fwd}
when $\overline{G}$ is Gorenstein with $a(\overline{G } ) = -1$.

\begin{proposition}
\label{proposition:LevelGaInvMinusOne}
Let $(R, \frakm)$ be a $d$-dimensional $F$-rational domain
of characteristic $p>0$ and $I$ an $\frakm$-primary ideal
such that $\overline{\Rees(I)}$ is Cohen-Macaulay.
Let $\overline{G} := \oplus_{n \geq 0}
\frac{\overline{I^n}}{\overline{I^{n+1}}}$.
Assume that
\[
\socle \left( \homology_{\overline{G}_+}^d(\overline{G})\right )
\subseteq
\left[ \homology_{\overline{G}_+}^d(\overline{G})\right ]_{-1}.
\]
Then $\overline{\Rees(I )}$ is $F$-rational.
\end{proposition}

A consequence of~\cite[Theorem~1.5]{HyryBlowupRingsRationalSings1999} is the
following:
Let $(A, \fraka)$ be an excellent local normal domain of characteristic zero
and $\frakb$ an $\fraka$-primary ideal of $A$.
Assume that the the normalized Rees algebra
$B := \oplus_{n \in \naturals}\overline{\frakb^n}$ is Cohen-Macaulay.
Then $B$ has rational singularities if and only if $\Proj B$ has rational
singularities. We wonder whether the following prime-characteristic
analogue is true. Assume, with notation as in the top of this section, that
$R$ is additionally $F$-rational;
if $\Proj \overline{\Rees(I )}$ is $F$-rational, is
$\overline{\Rees(I)}$ $F$-rational?
In this context, we have the following.

\begin{theorem}
\label{theorem:projimplies}
Let $(R,\mathfrak{m})$ be a complete $F$-finite normal domain of
characteristic $p \geq 7$ and with an infinite residue field, satisfying
one of the following:
\begin{enumerate}
\item
\label{theorem:projimplies:hyp}
$R$ is a hypersurface of dimension $d \geq 2$ and of multiplicity $2$;

\item
\label{theorem:projimplies:dim3}
$R$ is a three-dimensional Gorenstein $F$-rational ring.
\end{enumerate}
Suppose that $\Rees(\frakm)$ is Cohen-Macaulay and that
$\Proj \Rees(\frakm)$ is $F$-rational.
Then $\Rees(\frakm)$ is $F$-rational.
Additionally, in~\eqref{theorem:projimplies:hyp}, $R$ is $F$-rational.
\end{theorem}

This paper is organized as follows. In Section~\ref{section:prelims} we
collect relevant definitions and facts about Rees algebras, tight closure
and test ideals. Theorem~\ref{theorem:frationalityifftestideal} is proved
in Section~\ref{section:pfThmTestIdeal}. Proposition~\ref{proposition:gorG}
and other corollaries of Theorem~\ref{theorem:frationalityifftestideal}
are proved in Section~\ref{section:cor}.
Section~\ref{section:prelimsII} contains some background information for
the proof of Theorem~\ref{theorem:projimplies}, which is given in
Section~\ref{section:pfThmProjImplies}.

\subsection*{Acknowledgements}

We thank the referee for helpful comments.

\section{Preliminaries}
\label{section:prelims}
\numberwithin{equation}{subsection}

All the rings we consider in this article are excellent. The letter $p$
denotes a prime number. When used in the context of the Frobenius map and
singularities in prime characteristic, $q$ denotes an arbitrary power of
$p$.

We now collect
results on Rees algebras and singularities in prime characteristic that
will be required in the proofs.

\subsection{Local cohomology}
Let $R$ be a ring and $I = (r_1, \ldots, r_m)$ be an $R$-ideal. Then
the elements of $\homology_I^m(R)$ are the residue classes of the fractions
$\frac{a}{(r_1\cdots r_m)^l}$, $a \in R$ and $l \geq 1$,
modulo the $m$-boundaries in the
extended \tCech complex $\Cech^\bullet(r_1, \ldots, r_m)$.
When $r_1, \ldots, r_m$ is a regular sequence,
the class of $\frac{a}{(r_1\cdots r_m)^l}$ is zero if and only if $a \in
(r_1^l, \ldots, r_m^l)$;
see~\cite[Proof of Theorem 2.1, p.~104--105]{LipTeiPseudoRatlSing81}.

\subsection{Rees algebras}
Let $R$ be a ring.
A \define{filtration} of $R$ is a sequence
$\calI := (I_n)_{n \in \naturals}$ of $R$-ideals $I_n$ such that
\[
I_0 = R, I_{n+1} \subseteq I_n \;\text{and}\;
I_nI_m \subseteq I_{n+m} \;\text{for all}\; n, m \in \naturals.
\]
The \define{Rees algebra} of $\calI$ is the graded subring
\[
\Rees(\calI) := \bigoplus_{n \in \naturals} I_n t^n
\]
of $R[t]$ ($t$ being an indeterminate)
with $\deg r =0$ for each $r \in R$ and $\deg t = 1$.
Let $I$ be an ideal of $R$.
We write $\Rees(I)$ for the Rees algebra for the
$I$-adic filtration, i.e., the one with $I_n = I^n$ for each $n$.
We say that a filtration $\calI$ is \define{$I$-admissible}
if $I \subseteq I_1$ and $\Rees(\calI)$ is a finitely generated
$\Rees(I)$-module. By $\overline{\Rees(I)}$ we mean the integral
closure of $\Rees(I)$ in $R[t]$. It is the Rees algebra of the filtration
$(\overline{I^n})_{n \in \naturals}$, where for every $R$-ideal $J$,
$\overline{J}$ denotes its integral closure. Since $R$ is excellent,
the filtration $(\overline{I^n})_{n \in \naturals}$ is $I$-admissible.

A \define{reduction} of a filtration $\calI$ is an $R$-ideal $I \subseteq
I_1$ such that $\calI$ is $I$-admissible.
A \define{reduction} of an ideal $I$ is a reduction for the $I$-adic
filtration of $R$.
A reduction $I$ of
$\calI$  is \define{minimal} if it is minimal with respect to inclusion.
If $R$ is a local ring with maximal ideal $\frakm$ and $R/\frakm$ is
infinite, then every minimal reduction of $\calI$ is minimally generated by
$\dim(\Rees(\calI) \otimes_R R/\frakm)$ elements. In particular, if $I_1$
is $\frakm$-primary, then
every reduction of $\calI$ is minimally generated by $\dim R$ elements.
For a reduction $I$ of $\calI$,
the \define{reduction number} $r_I(\calI)$ is the minimum of the integers
$r$ such that $I_{n+1} = II_n$ for all $n \geq r$.
When $\calI$ is the $I$-adic filtration (for some $R$-ideal $I$)
and $J$ is a reduction of $I$, we write $r_J(I) = r_J(\calI)$.

The \emph{associated graded ring} of a filtration $\calI = (I_n)$ is
$G(\calI) := \bigoplus_{n\geq 0}\frac{I_n}{I_{n+1}}$.
When $\calI$ is the $I$-adic filtration, we write $G(I) = G(\calI)$;
when $\calI = (\overline{I^n})_{n\geq 0 }$, we write
$\overline G(I) = G(\calI)$.

We now collect some statements about the local cohomology modules of Rees
algebras. Recall that for a noetherian graded ring
$S := \oplus_{n \geq 0 } S_n$ with
$S_0$ local, the \define{$a$-invariant} $a(S)$ of $S$ is
\[
\max \{ j \mid \left[\homology^{\dim S}_\frakn(S)\right]_j \neq 0\}
\]
where $\frakn$ is the unique homogeneous maximal ideal of $S$.
If $S$ is a $\ints$-graded ring with unique homogeneous maximal ideal
$\frakn$ and $M$ is a graded $S$-module, the \define{socle} of $M$,
denoted $\socle M$ is $(0 :_M \frakn)$.

\begin{theorem}
\label{theorem:cmproperties}
Let $(R,\mathfrak{m})$ be a $d$-dimensional Cohen-Macaulay local ring
and $I$ an $\mathfrak{m}$-primary ideal.
Let $\calI$ be an $I$-admissible filtration.
\begin{enumerate}%

\item
\label{theorem:cmproperties:ainv}
\cite[Lemma~3.3]{GotoNishidaFiltrGorRees1994}
$a(\Rees(\calI)) = -1$.

\item
\label{theorem:cmproperties:cm}
\cite[Corollary 1.2]{GotoNishidaFiltrGorRees1994}
$\Rees(\calI)$ is Cohen-Macaulay if and
only if $G(\calI)$ is Cohen-Macaulay and
$a(G(\calI))<0$.
\end{enumerate}
\end{theorem}

In some cases, there is the following relation between the $a$-invariant
and the reduction number, which can be proved similar
to~\cite[Proposition~3.2]{TrungRedExpDegBd87}.
(The second assertion follows from the first and 
Theorem~\ref{theorem:cmproperties}~\eqref{theorem:cmproperties:cm}.)

\begin{proposition}
\label{proposition:aInvRedNo}
With notation as in Theorem~\ref{theorem:cmproperties}, assume that
$G(\calI)$ is Cohen-Macaulay.
Let $I$ be a minimal reduction of $\calI$.
Then $r_I(\calI) = a(G(\calI)) + d$.
In particular, if $\Rees(\calI)$ is Cohen-Macaulay,  
then $r_I(\calI) \leq d-1$.
\end{proposition}

Write $\Rees'(\calI)$ for the extended Rees algebra
$\oplus_{n \in \ints} I_nt^n$ (where $I_n := R$ when $n \leq 0$).
W.~Heinzer, M.-K.~Kim and
B.~Ulrich~\cite[Theorem~6.1]{HeinzerKimUlrichCMGor2011}
gave a description of $\omega_{\Rees'(\calI)}$. Together
with~\cite[Discussion~3.1]{KoleyKumminiFrationality2017}, we get the
following proposition.

\begin{proposition}
\label{proposition:descrOmegaR}
Let $(R,\mathfrak{m})$ be a Gorenstein local ring, of dimension
$ d\geq 2$.
Let $I$ be an ideal generated by a system of parameters and $\calI$ be an
$I$-admissible filtration.
Write $r = r_I(\calI)$.
If $\Rees(\calI)$ is Cohen-Macaulay, then
$\left[\omega_{\Rees(\calI)}\right]_n=I^n:I_{d-1}$
for each $n \geq 1$.
If, further, $G(\calI)$ is Gorenstein, then
$\left[\omega_{\Rees(\calI)}\right]_n=I_{n+a+1}$
where $a=a(G(\calI))$.
\end{proposition}

\begin{proof}
By~\cite[Theorem~6.1(1)]{HeinzerKimUlrichCMGor2011},
\[
\left[\omega_{\Rees'(\calI)}\right]_n=I^{n-(d-1)+r}:I_{r}
\]
for each $n$.
Since $\Rees(\calI)$ is Cohen-Macaulay, 
$r \leq d-1$
(Proposition~\ref{proposition:aInvRedNo}), so
$I^{n}:I_{d-1}= I^{n}: I^{(d-1)-r}I_{r} = I^{n-(d-1)+r}:I_{r}$.
By~\cite[Discussion~3.1]{KoleyKumminiFrationality2017},
for each $n < 0$,
\[
\left[\homology_{\frakM}^{d+1}(\Rees(\calI))\right]_n
=
\left[\homology_{\frakM'}^{d+1}(\Rees'(\calI))\right]_n
\]
where $\frakM$ and $\frakM'$ are the maximal homogeneous ideals of
$\Rees(\calI)$ and $\Rees'(\calI)$ respectively.
(In each negative degree, the middle vertical map
in the commutative diagram
in~\cite[Discussion~3.1]{KoleyKumminiFrationality2017} is the identity
map, and, therefore, so is the map $\gamma$.)
Applying the graded Matlis duality functor $\Hom_R(-, E)$
(where $E$ is the injective hull of the $R$-module $R/\frakm$)
we see that for each $n > 0$,
\[
\left[\omega_{\Rees(\calI)}\right]_n =
\left[\omega_{\Rees'(\calI)}\right]_n =
I^{n}:I_{d-1}.
\]
Assume further that $G(\calI)$ is Gorenstein.
Then $\omega_{G(\calI)} = G(\calI)(a)$.
Since $G(\calI)=\frac{\Rees'(\calI )}{t^{-1}\Rees'(\calI )}$,
it follows that $\omega_{\Rees'(\calI ) } = \Rees'(\calI )(a+1)$.
\end{proof}

\subsection{Prime characteristic notions}

Let $R$ be a ring of prime characteristic $p>0$.
For an ideal $I\subseteq R$, $I^{[q]}$ denotes the
ideal generated by $q$-th power of elements of $I$.
Let $F:R\to R$ be the \emph{Frobenius map} $r \mapsto r^p$.
The ring $R$ viewed as an $R$-module via $F^e:R\to R$ is denoted as
${}^e\!R$.
Let $M$ be an $R$-module and write $F^e(M) = {}^e\!R\otimes_RM$.
For $x \in M$, its image in $F^e(M)$ under the
natural map $M\to {{}^e\!R}\otimes_R M$ is written as $x^{p^e}$.
For a submodule $N\subseteq M$, the image of $F^e(N)$
in $F^e(M)$ is denoted by $N^{[p^e]}_M$.
Write
\[
R^{0} = R \minus \bigcup_{\frakp \in \Min R }\frakp.
\]

Let $\fraka$ be an $R$-ideal and $N \subseteq M$ $R$-modules.
The \define{$\mathfrak{a}$-tight closure} of $N$ in $M$
(introduced by N.~Hara and
K.-i.~Yoshida~\cite[Definition~1.1]{HaraYoshidaGenTightClMultIdeals2003})
is
\[
N^{*\mathfrak{a}}_M :=
\{z\in M \mid
\text{there exists}\;
c\in R^{0} \;\text{such that}\; c\mathfrak{a}^qz^q\subseteq N^{[q]}_M
\;\text{for all}\; q\gg 1 \}.
\]
When $\fraka = R$, this is the same as the tight closure defined by
M.~Hochster and C.~Huneke\cite{HochsterHunekeTCInvThyBSThm90}; we then
write $N^{*}_M$ for $N^{*R}_M$.
When $M=R$ and $N=I$ an $R$-ideal, we
write $I^{*\mathfrak{a}}$ and $I^*$ respectively.
We say that $N$ is \define{tightly closed} in $M$ if
$N^{*}_M=N$.

\begin{definition}
We say that a local ring $R$ is \define{$F$-rational} if an ideal generated
by a system of parameters is tightly closed.
We say that a ring $R$ is \define{$F$-rational} if $R_\frakm$ is
$F$-rational for each maximal ideal $\frakm$ of $R$.
Say that a scheme $X$ is \define{$F$-rational} if it has an open cover by
the spectra of $F$-rational rings.
\end{definition}

We have used~\cite[Theorem~(4.2)]{HochsterHunekeFregSmBsch94} to give the
above definition of $F$-rationality.
Moreover, the same theorem also shows that all $F$-rational rings are
normal, and those that we consider in this paper are Cohen-Macaulay.
A Gorenstein local $F$-rational ring is weakly $F$-regular, i.e., all
ideals are tightly closed.
It is known that a Cohen-Macaulay ring $R$ is $F$-rational if and only if
$0^{*}_{H_{\mathfrak{m}}^d(R)} = 0$~\cite[Theorem~2.6]{SmithFrational1997}.
Suppose that $\mathscr{R}(\calI)$  is Cohen-Macaulay.
Then by a similar argument, $\mathscr{R}(\calI)$ is $F$-rational if and
only if $0^{*}_{H_{\mathfrak{M}}^{d+1}(\mathscr{R}(\calI))}=0$.

\begin{remarkbox}
\label{remarkbox:reesnormal}
A Cohen-Macaulay domain of dimension at least two that is regular in
codimension one is normal.
Suppose that $\Rees(\calI)$ is Cohen-Macaulay and that $\Proj
\Rees(\calI)$ is $F$-rational.
As discussed in \cite[Lemma~3.5]{KoleyKumminiFrationality2017},
$\Spec\Rees(\calI)\setminus \{\mathfrak{M}\}$
(where $\frakM$ is the homogeneous maximal ideal of $\Rees(\calI)$)
is $F$-rational.
Therefore $\Rees(\calI)$ is normal.
\end{remarkbox}

An element $c\in R^{0}$ is called \define{$\mathfrak{a}$-test element} if for
every ideal $I$ of $R$ and for all $z\in R$, $z\in I^{*\mathfrak{a}}$ if
and only if $cz^q\mathfrak{a}^q\subseteq I^{[q]}$ for all $q\geq 1$.
An $R$-test element is a \define{test element}.
An element $c\in R^{0}$ is called a \define{parameter test element} if
for every $R$-ideal $I$ generated by a system of parameters and for all
$z\in R$, $z\in I^{*}$ if and only if
$cz^q\subseteq I^{[q]}$ for all $q\geq 1$.
The \define{test ideal} $\tau (\mathfrak{a})$ is $\bigcap_{I\subseteq
R}(I:I^{*\mathfrak{a}})$.
If $(R,\mathfrak{m})$ is a $d$-dimension normal and Gorenstein local ring,
then
\begin{equation}
\label{equation:tauGor}
\tau(\mathfrak{a})=\operatorname{Ann}_R(0^{*\mathfrak{a}}_E)=\bigcap_{t\geq
1}(x_1^t,...,x_d^t):_R(x_1^t,...,x_d^t)^{*\mathfrak{a}}.
\end{equation}
where $x_1, \ldots, x_d$ is a system of parameters and $E$ is the
\emph{injective hull} of $R/\mathfrak{m}$.
Recall that
$E \simeq \homology^d_{\mathfrak{m}}(R)
\simeq \lim\limits_{\rightarrow  t} \frac{R}{(x_1^t,...,x_d^t)}$.

We summarize some results
from~\cite[Proposition~1.11 and
Theorem~2.1]{HaraYoshidaGenTightClMultIdeals2003}
that are relevant in this paper.

\begin{theorem}
\label{theorem:testIdealProps}
Assume the above notation.
Let $\mathfrak{a},\mathfrak{b}$ be $R$-ideals.
Then:

\begin{enumerate}

\item
\label{theorem:testIdealProps:prod}
$\fraka \tau(\mathfrak{b}) \subseteq \tau(\mathfrak{a}\frakb)$.

\item
If $\mathfrak{b}\subseteq\mathfrak{a}$, then $\tau(\mathfrak{b})\subseteq
\tau(\mathfrak{a})$.
Moreover, if $\mathfrak{a}\cap R^0\neq \emptyset$ and $\mathfrak{b}$ is a
reduction of $\mathfrak{a}$, then
$\tau(\mathfrak{b})= \tau(\mathfrak{a})$.

\item If $R$ is a weakly $F$-regular ring, then $\mathfrak{a}\subseteq
\tau(\mathfrak{a})$.

\end{enumerate}
\end{theorem}

\section{Proof of Theorem~\protect{\ref{theorem:frationalityifftestideal}}}
\label{section:pfThmTestIdeal}
\numberwithin{equation}{section}

\begin{setup}
\label{setup:main}
Let $(R,\mathfrak{m})$ be an excellent, Cohen-Macaulay, normal local domain
of dimension $d\geq 2$ and of characteristic $p>0$.
Let $I=(f_1,\cdots,f_d)$ be a parameter ideal, that is $f_1,\cdots ,f_d$ is
a system of parameters (and hence a regular sequence).
Let $f=f_1\cdots f_d$.
Set $I^{[l]}:=(f_1^l,\cdots,f_d^l)$ for $l\geq 1$.
Denote the Rees algebra $R[It]$ (respectively $\overline{R[It]}$) by
$\Rees$ (respectively $\overline{\Rees}$).
$\mathfrak{M}$ (respectively $\overline{\mathfrak{M}}$) to be the unique
homogeneous maximal ideal of $\Rees$ (respectively $\overline{\Rees}$).
$\overline{G}$ is the associated graded ring $\bigoplus_{n\geq
0}\frac{\overline{I^n}}{\overline{I^{n+1}}}$.
\end{setup}

The following proposition is similar to those
in~\cite{ValabregaVallaFormRings1978}. We could not find a proof of the
statement as we need it below, so we give a proof.

\begin{proposition}
\label{proposition:VVlike}
Assume Setup~\ref{setup:main} and that $\overline{G}$ is Cohen-Macaulay.
Then for all $k \geq 0$ and $l \geq 0$,
$\overline{I^{k+l}} \cap I^{[l]} = \overline{I^{k}} I^{[l]}$.
\end{proposition}

\begin{proof}
For brevity of notation, we write $x_i$ for the image of $f_i$ in
$\overline{G}_1$.
Let $a_i, 1 \leq i \leq d$ be elements of $R$ such that
$\sum_{i}^d a_i f_i^l \in \overline{I^{k+l}}$.
Let $k_1$ be such that $(a_1, \ldots, a_d ) \in
\overline{I^{k_1}} \minus \overline{I^{k_1+1}}$.
Let $\bar{a_i}, 1 \leq i \leq d$ be the images of the $a_i$ in
$\overline{G}_{k_1}$.
(We write $\bar a_i$ instead of $a_i^*$ since not all the $a_i$ might be in
$\overline{I^{k_1}} \minus \overline{I^{k_1+1}}$.)
We may assume that $k_1 < k$.

Let $(K_\bullet, \partial_\bullet)$
(respectively, $(K'_\bullet, \delta_\bullet)$) be the
Koszul complex on the $f_i^l$ on $R$  (respectively, the $x_i^l$ on
$\overline{G}$).
Note that $x_1, \ldots, x_d$ is a regular sequence on $\overline{G}$.
Therefore
\[
\begin{bmatrix}
\bar a_1 \\ \bar a_2 \\  \vdots \\ \bar a_d
\end{bmatrix}
\in \ker \delta_1 = \image \delta_2.
\]
Hence there exist $b_1, \ldots, b_{\binom{d}{2}} \in \overline{I^{k_1-l}}$
such that
\[
\begin{bmatrix}
\bar a_1 \\ \bar a_2 \\  \vdots \\ \bar a_d
\end{bmatrix}
=
\delta_2
\left(
\begin{bmatrix}
\bar b_1 \\ \bar b_2 \\  \vdots \\ \bar b_{\binom{d}{2}}
\end{bmatrix}
\right)
\]
Define $a'_i, 1 \leq i \leq d$ by
\[
\begin{bmatrix}
a'_1 \\  a'_2 \\  \vdots \\  a'_d
\end{bmatrix}
=
\begin{bmatrix}
 a_1 \\  a_2 \\  \vdots \\  a_d
\end{bmatrix}
-
\partial_2
\left(
\begin{bmatrix}
 b_1 \\  b_2 \\  \vdots \\  b_{\binom{d}{2}}
\end{bmatrix}
\right)
\]
Note that $a'_i \in \overline{I^{k_1+1}}$ for all $i$.
Moreover,
\[
\sum_i a'_i f_i^l =
\partial_1
\left(
\begin{bmatrix}
a'_1 \\  a'_2 \\  \vdots \\  a'_d
\end{bmatrix}
 \right)
=
\partial_1
\left(
\begin{bmatrix}
a_1 \\  a_2 \\  \vdots \\  a_d
\end{bmatrix}
 \right)
 =
\sum_i a_i f_i^l.
\]
Repeating this argument, we can find $\tilde a_i \in \overline{I^{k}}$ such
that $\sum_i a_i f_i^l= \sum_i \tilde a_i f_i^l$.
\end{proof}

\begin{discussion}
\label{discussion:jlstar}
Assume Setup~\ref{setup:main}.
Then $I^{[l]}:_R I^{n}=I^{dl-n-d+1}+I^{[l]}$ for all $l,n\geq 1$;
see, e.g.,~\cite[Lemma 2.11]{HaraYoshidaGenTightClMultIdeals2003}.
Using this, Hara and Yoshida~\cite[Theorem
2.7]{HaraYoshidaGenTightClMultIdeals2003}
showed that
\begin{equation}
\label{equation:jlstardefnHY}
I^{[l]^{*I^n}} =\{z\in R: \;\text{there exists a nonzero}\; c\in R
\;\text{such that}\; cz^q\in I^{(dl-n)q}+I^{[lq]} \;\text{for
all}\; q\gg 1\}.
\end{equation}

Let $z \in R$, $c \in R \minus \{0\}$ and $q$ be such that
$cz^q\in \overline{I^{(dl-n)q}}+I^{[lq]}$.
Let $r = r_I((\overline{I^n})_{n \geq 0 })$.
Then for all $n' \geq 0$, $\overline{I^{n'+r}} \subseteq I^{n'}$.
Let $c_1 \in I^r \minus \{0 \}$.
Then
\[
c_1cz^q \in I^r\overline{I^{(dl-n)q}} + I^{[lq ]} \subset
I^{(dl-n)q}+  I^{[lq]}.
\]
Therefore,
\begin{equation}
\label{equation:jlstardefn}
I^{[l]^{*I^n}} =\{z\in R: \;\text{there exists a nonzero}\; c\in R
\;\text{such that}\; cz^q\in \overline{I^{(dl-n)q}}+I^{[lq]} \;\text{for
all}\; q\gg 1\},
\end{equation}
where $l,n\geq 1$ such that $dl-n\geq 1$.
In particular,
\begin{equation}
\label{equation:intCljlStar}
\overline{I^{dl-n}}+I^{[l]}\subseteq I^{[l]^{*I^n}}
\end{equation}
for all $l,n\geq 1$ such that $dl-n\geq 1$.
\end{discussion}

\begin{discussionbox}
\label{discussionbox:descrLC}
We now make some observations about
$\homology^d_{\overline{\Rees}_{+}} (\overline{\Rees})$ that are analogous
to those in
\cite[Subsection~1.3, Lemma~2.8,
Corollary~2.9]{HaraWatanabeYoshidaFrationality2002} about
$\homology^d_{\Rees_{+}} (\Rees)$.
Assume that $\overline{\Rees }$ is Cohen-Macaulay.
Let $n \in \ints$.
Then
\[
\left[ \homology^d_{\overline{\Rees}_{+}} (\overline{\Rees}) \right]_{n} =
\left \{
\left[\frac{a}{f^l}t^{n}\right] : l\geq 0, dl+n\geq 0, a \in
\overline{I^{dl+n}} \right \}.
\]
Assume now that $n \geq 1$ and consider the exact sequence
\begin{equation}
\label{equation:locohExactSeq}
\begin{tikzcd}
0 \arrow[r] &
\left[ \homology^d_{\overline{\Rees}_{+}} (\overline{\Rees})
\right]_{-n} \arrow[r, "\phi_{-n}"]
&
\homology^d_\mathfrak{m} (R)t^{-n} \arrow[r]
&
\left[ \homology^{d+1}_{\mathfrak{M}} (\overline{\Rees}) \right]_{-n}
\arrow[r]
& 0 .
\end{tikzcd}
\end{equation}
(The existence of such an exact sequence can be proved in a way similar
to~\cite[Lemma~2.7]{HaraWatanabeYoshidaFrationality2002}.)
The map $\phi_{-n}$ in~\eqref{equation:locohExactSeq} is
$\left[\frac{a}{f^l}t^{-n}\right] \mapsto
\left[\frac{a}{f^l}\right]t^{-n}$.
\end{discussionbox}

\begin{proposition}
\label{proposition:imagePhi}
Assume Setup~\ref{setup:main} and that $\overline{\Rees }$ is Cohen-Macaulay.
Suppose that $l,n$ are positive integers such that $dl-n\geq 1$.
Then for all $a\in R$,
$\left[\frac{a}{f^l}\right]t^{-n}\in \image (\phi_{-n})$
if and only if $a\in \overline{I^{dl-n}}+I^{[l]}$.
\end{proposition}

\begin{proof}
`If' follows from Discussion~\ref{discussionbox:descrLC}.
For `only if',
let $a \in R$ be such that
\[
\left[\frac{a}{f^l}\right]t^{-n}\in \image(\phi_{-n}).
\]
Then there exist $m \geq 0$ and $b \in \overline{I^{dm-n}}$ such that
$dm-n \geq 0$ and
\[
\left[\frac{a}{f^l}\right]t^{-n}
=
\left[\frac{b}{f^m}\right]t^{-n}
=
\phi_{-n}\left (\left[\frac{b}{f^m}t^{-n}\right]\right).
\]
We see that
\[
\frac{a}{f^l}-\frac{b}{f^m}
\]
is a boundary in $\Cech^\bullet(f_1, \ldots, f_d; R)$.
Without loss of generality, $m \geq l$.
Therefore there exists $N \geq m$ such that
\[
\frac{a}{f^l}-\frac{b}{f^m} = \frac{\sum_i a_i f_i^N }{f^N }.
\]
Hence $f^{N-m}(af^{m-l }-b ) \in I^{[N]}$
so
$af^{m-l} \in \overline{I^{dm-n}} + I^{[m]}$.
From this we want to conclude that
\[
a \in \overline{I^{dl-n}} + I^{[l]}.
\]
We may apply induction on $m-l$ and assume that $m=l+1$.

Let $k$ be such that $a \in \overline{I^k} \minus \overline{I^{k+1}}$.
We may assume that $k < dl-n$.
Further applying induction on $dl-n-k$, we may further assume that
\begin{equation}
\label{equation:largerk}
((\overline{I^{d{l+1}-n}} + I^{[{l+1}]}) : f) \cap \overline{I^{k'}}
 \subseteq
\overline{I^{dl-n}} + I^{[l]}.
\end{equation}
for all $k'$ with $k < k' \leq dl-n$.

Write
\begin{equation}
\label{equation:afisbplussth}
af = b + \sum_{i=1}^d a_i f_i^{l+1}
\end{equation}
where $b \in \overline{I^{d(l+1)-n}}$ and $a_i \in R$ for all $1 \leq i
\leq d$.
By Proposition~\ref{proposition:VVlike}, we may
assume that $a_i \in \overline{I^{k+d-l-1}}$ for all $1 \leq i \leq d$.

Write $x_i$ for the image of $f_i$ in $\overline{G}_1$
and $\bfx = x_1\cdots x_d$.
Since $\overline{G}$ is Cohen-Macaulay, the $x_i$ form a regular sequence
on it.
Then we get a map between the Koszul co-complexes
$K^\bullet(x_1^l, \ldots, x_d^l; \overline{G}) \to
K^\bullet(x_1^{l+1}, \ldots, x_d^{l+1}; \overline{G})$
\[
\begin{tikzcd}
\cdots \arrow[r] &
\overline{G}((d-2)l)^{\binom{d}{2}} \ar[r] \ar[d]&
\overline{G}((d-1)l)^d \ar[r, "\partial^{(l)}"] \ar[d]&
\overline{G}(dl) \ar[r] \ar[d, "\cdot \bfx"]& 0
\\
\cdots \arrow[r] &
\overline{G}((d-2)(l+1))^{\binom{d}{2}} \ar[r] &
\overline{G}((d-1)(l+1))^d \ar[r, "\partial^{(l+1)}"] &
\overline{G}(d(l+1)) \ar[r] & 0
\end{tikzcd}
\]
Note that $\deg a^*\bfx = k+d < d(l+1)-n = \deg b^*$, so
from~\eqref{equation:afisbplussth}
we see that $a^*\bfx \in \image \partial^{(l+1)}$, i.e. it is boundary.
Hence it gives the zero element in
$\homology_{\overline{G}_+}^d(\overline{G})$.
Since $\overline{G}$ is Cohen-Macaulay, all the maps in the directed system
\[
\homology^d(K^\bullet(x_1^l, \ldots, x_d^l; \overline{G}))
\to
\homology^d(K^\bullet(x_1^{l+1}, \ldots, x_d^{l+1}; \overline{G}))
\]
are injective, so $a^*$ too is boundary, i.e., $a^* \in \image
\partial^{(l)}$. I.e., there exist $\alpha_1, \ldots, \alpha_d \in R$ such
that
\[
a' := a - \sum_{i} \alpha_i f_i^l \in \overline{I^{k+1}}
\]
Note that $a'f - af \in I^{[l+1]}$
so by~\eqref{equation:largerk},
$a' \in \overline{I^{dl-n}} + I^{[l]}$.
Hence
$a \in \overline{I^{dl-n}} + I^{[l]}$.
\end{proof}

\begin{discussionbox}
\label{discussion:hwydiscussiononexactness}
Since $R$ is an excellent normal domain, its singular locus is a proper
closed subset of $\Spec R$, so there exists
$c \in R^0$ such that $R_c$ is regular, and, \textit{a fortiori},
$F$-rational. Then
$\left({\Rees}\right)_c \simeq \left(\overline{\Rees}\right)_c \simeq R_c[t]$
is $F$-rational; use~\cite[Proposition~1.2]{VelezOpennessFrationalLoci1995}.
Hence by~\cite[Theorem~3.9]{VelezOpennessFrationalLoci1995}, we see that
there exists $N$ such that $c^N$ is a parameter test element for
$R$, ${\Rees}$ and $\overline{\Rees}$.
In particular, $I$ contains parameter test elements for
$R$, ${\Rees}$ and $\overline{\Rees}$.

Let $c\in I$ a parameter test element for $R$ and $\overline{\Rees}$.
Then for all $n\geq 1$, we have, from~\eqref{equation:locohExactSeq},
the following commutative diagram:
\begin{equation}
\label{equation:locohExactSeqFrob}
\begin{tikzcd}
0 \arrow[r] & \left[ \homology^d_{\overline{\Rees}_{+}} (\overline{\Rees})
\right]_{-n} \arrow[r, "\phi_{-n}"] \arrow[d,"cF^e"]
& \homology^d_\mathfrak{m} (R)t^{-n} \arrow[r] \arrow[d, "cF^e"]
&\left[ \homology^{d+1}_{\mathfrak{M}} (\overline{\Rees}) \right]_{-n}
\arrow[r] \arrow[d, "cF^e"]  & 0
\\
0 \arrow[r] & \left[ \homology^d_{\overline{\Rees}_{+}} (\overline{\Rees})
\right]_{-nq} \arrow[r, "\phi_{-nq}"]
& \homology^d_{\mathfrak{m}} (R)t^{-nq} \arrow[r] & \left[
\homology^{d+1}_{\mathfrak{m}} (\overline{\Rees}) \right]_{-nq} \arrow[r]
& 0 .
\end{tikzcd}
\qedhere
\end{equation}
\end{discussionbox}

\begin{lemma}
\label{lemma:jlstarnondecreasing}
Let $(R,\mathfrak{m})$ and $I$ be as in the Setup~\ref{setup:main}.
Then for each $n\geq 1$ and $l\geq 1$,
$I^{[l+1]}:I^{[l+1]^{*I^n}} \subseteq I^{[l]}:I^{[l]^{*I^n}}$.
\end{lemma}

\begin{proof}
Let $n \geq 1$ and $l \geq 1$.
Note that $I^{[l+1]}:f=I^{[l]}$ since
$f_1,f_2,\cdots, f_d$ is a regular sequence.
Hence we have the following sequence of inclusions:
\begin{align*}
fI^{[l]^{*I^n}} & \subseteq I^{[l+1]^{*I^n}};\\
fI^{[l]^{*I^n}}(I^{[l+1]}:I^{[l+1]^{*I^n}}) & \subseteq
I^{[l+1]^{*I^n}}(I^{[l+1]}:I^{[l+1]^{*I^n}}) \subseteq I^{[l+1]};\\
I^{[l]^{*I^n}}(I^{[l+1]}:I^{[l+1]^{*I^n}}) & \subseteq I^{[l+1]}:f=I^{[l]}.
\end{align*}
Therefore
$I^{[l+1]}:I^{[l+1]^{*I^n}} \subseteq I^{[l]}:I^{[l]^{*I^n}}$.
\end{proof}

\begin{lemma}
\label{lemma:vanishingzerostar}
Let $(R,\mathfrak{m})$ and $I$ as in the Setup~\ref{setup:main}. Suppose
that $\overline{\Rees}$ is Cohen-Macaulay. Fix $n\geq 1$. Consider
the following statements:
\begin{enumerate}%

\item
\label{lemma:vanishingzerostar:zerostar}
$\left[ 0^{*}_{\homology^{d+1}_{\overline{\mathfrak{M}}} (\overline{\Rees})}
\right]_{-n} = 0$.

\item
\label{lemma:vanishingzerostar:jlstar}
$I^{[l]^{*I^n}} =
\overline{ I^{dl-n}}+I^{[l]}$ for all $l$ such that $dl-n\geq  1$.

\item
\label{lemma:vanishingzerostar:jlstarlargel}
$I^{[l]^{*I^n}}=\overline{I^{dl-n}}+I^{[l]}$ for all $l\gg 0$.

\item
\label{lemma:vanishingzerostar:tau}
$\tau (I^n)=I^n:\overline{ I^{d-1}}$.
\end{enumerate}

Then \eqref{lemma:vanishingzerostar:zerostar} $\iff$
\eqref{lemma:vanishingzerostar:jlstar} $\iff$
\eqref{lemma:vanishingzerostar:jlstarlargel}.
If, additionally, $R$ is Gorenstein, then
\eqref{lemma:vanishingzerostar:jlstarlargel} $\iff$
\eqref{lemma:vanishingzerostar:tau}.
\end{lemma}

\begin{proof}
\eqref{lemma:vanishingzerostar:zerostar} $\implies$
\eqref{lemma:vanishingzerostar:jlstar}:
Let $l$ be such that $dl-n\geq 1$.
In view of~\eqref{equation:intCljlStar}, we need to show that
$I^{[l]^{*I^n}} \subseteq \overline{ I^{dl-n}}+I^{[l]}$.
Let $a\in I^{[l]^{*I^n}}$.
Then there exists a non-zero $c \in R$ such that
$ca^q\in I^{(dl-n)q}+I^{[lq]} \subseteq
\overline{ I^{(dl-n)q}}+I^{[lq]}$
for all $q\gg 1$, by~\eqref{equation:jlstardefnHY}.
Hence by Proposition~\ref{proposition:imagePhi},
\begin{equation}
\label{equation:imagePhinq}
cF^e\left(\left[\frac{a}{f^l}\right]t^{-n}\right) =
\left[\frac{ca^{p^e}}{f^{l{p^e}}}\right]t^{-n{p^e}}
\in \image \phi_{-n{p^e}}
\end{equation}
for all $e \gg 0$.
Multiplying $c$ by a parameter test element in $I$ for
$R$ and $\overline{\Rees}$
(Discussion~\ref{discussion:hwydiscussiononexactness}),
we may assume that $c$ is a parameter test element.
Consider the element
$\left[\frac{a}{f^l}\right]t^{-n}\in \homology^d_\mathfrak{m} (R)t^{-n}$ and
its image $\xi$ in
$[\homology^{d+1}_{\overline{\mathfrak{M}}} (\overline{\Rees})]_{-n}$.
From~\eqref{equation:imagePhinq}
and the commutative diagram~\eqref{equation:locohExactSeqFrob},
we see that $cF^e(\xi)=0$ for all $e\gg 0$.
Hence
$\xi \in \left[ 0^{*}_{\homology^{d+1}_{\overline{\mathfrak{M}}}
(\overline{\Rees})} \right]_{-n}$ and thus
$\xi =0$.
So $[\frac{a}{f^l}]t^{-n}\in \image(\phi_{-n})$.
Therefore
$a\in \overline{I^{(dl-n)}}+I^{[l]}$.

\eqref{lemma:vanishingzerostar:jlstar} $\implies$
\eqref{lemma:vanishingzerostar:jlstarlargel}:
Immediate.

\eqref{lemma:vanishingzerostar:jlstarlargel} $\implies$
\eqref{lemma:vanishingzerostar:zerostar}:
Let $\xi \in
\left[0^{*}_{\homology^{d+1}_{\overline{\mathfrak{M}}} (\overline{\Rees})}
\right]_{-n}$.
Let $a \in R$ and $l \in \naturals$ be such that $\xi$ is the image of
the element
\[
\left[\frac{a}{f^l}\right]t^{-n} \in
\homology^d_\mathfrak{m} (R) t^{-n}.
\]
Since
$\left[\frac{a}{f^l}\right] = \left[\frac{af^{l^{'}}}{f^{l+l^{'}}}\right]$
in $\homology^d_{\mathfrak{m}}(R)$ for all $l^{'}\geq 0$,
we may assume that $l$ is sufficiently large. Now $cF^e(\xi)=0$ for all
$e\gg 0$. From the commutative diagram~\eqref{equation:locohExactSeqFrob},
we see that
\[
\left[\frac{ca^q}{f^{lq}}\right]t^{-nq} \in \image(\phi_{-nq})
\]
so by Proposition~\ref{proposition:imagePhi}
and~\eqref{equation:jlstardefn},
$a\in I^{[l]^{*I^n}}$.
However, $I^{[l]^{*I^n}} = \overline{I^{dl-n}}+I^{[l]}$, by hypothesis.
Therefore
$\left[\frac{a}{f^l}\right]\in \image (\phi_{-n})$ and so $\xi =0$.

We now prove
\eqref{lemma:vanishingzerostar:jlstarlargel} $\iff$
\eqref{lemma:vanishingzerostar:tau} assuming that $R$ is Gorenstein.
Note that
\begin{equation}
\label{equation:colon}
I^{[l]}:(\overline{I^{dl-n}}+I^{[l]})
=I^{[l]}:\overline{I^{dl-n}}
=(I^{[l]}:I^{dl-n-d+1}):\overline{I^{d-1}}
=(I^{[l]}+I^n):\overline{I^{d-1}}
=I^n:\overline{I^{d-1}}.
\end{equation}
for all $l \geq n$.
(We have used the fact that
$\overline{I^{dl-n}} = I^{dl-n-d+1}\overline{I^{d-1}}$, 
which, in turn, follows from the fact that $r_I(\overline{I}) \leq d-1$, by
Proposition~\ref{proposition:aInvRedNo}.)
It immediately follows that if we
assume~\eqref{lemma:vanishingzerostar:jlstarlargel}, then,
using~\eqref{equation:tauGor} and Lemma~\ref{lemma:jlstarnondecreasing},
\[
\tau (I^n) = \bigcap_{l\geq 1} I^{[l]}:I^{[l]^{*I^n}}
= I^n:\overline{I^{d-1}}.
\]

Conversely assume~\eqref{lemma:vanishingzerostar:tau}.
In view of~\eqref{equation:intCljlStar}, we need to show that
$I^{[l]^{*I^n}} \subseteq \overline{I^{dl-n}}+I^{[l]}$
for all $l\gg 0$.
Or, equivalently, using~\eqref{equation:colon}, that
$I^{[l]} : I^{[l]^{*I^n}} \supseteq I^n:\overline{I^{d-1}}$
for all $l\gg 0$.
This follows from~\eqref{equation:tauGor}.
\end{proof}

\begin{remark}
\label{remark:tauInOmegan}
Adopt the notation and hypotheses of
Theorem~\ref{theorem:frationalityifftestideal}.
By~\cite[Theorem~2.1]{HaraYoshidaGenTightClMultIdeals2003},
$\overline{I^{d-1}} \tau(I^n) \subseteq \tau (I^{n+d-1}) \subseteq I^n$.
Hence $\tau(I^n) \subseteq
\left[\omega_{\overline{\Rees}}\right]_n$, using
Proposition~\ref{proposition:descrOmegaR}.
Thus, together
with~\cite[Theorem~5.1]{HaraYoshidaGenTightClMultIdeals2003},
Theorem~\ref{theorem:frationalityifftestideal}
implies that
$\overline{\Rees}$ is $F$-rational if and only if
$\tau(I^n) = \left[\omega_{\overline{\Rees}}\right]_n$ for all
$n\geq 1$.
\end{remark}

\begin{proof}%
[Proof of Theorem~\protect{\ref{theorem:frationalityifftestideal}}]
We first show that
$\tau(I^n)=I^n:\overline{ I^{d-1}}$ for each $n\geq 1$.
By hypothesis this holds for all $1 \leq n \leq d-1$.
Now assume that $n \geq d$, by induction, that
$\tau(I^n)=I^n:\overline{ I^{d-1}}$, which, in turn, equals
$\left[\omega_{\overline{\Rees}}\right]_n$.
Therefore
\[ \tau(I^{n+1})\subseteq \left[\omega_{\overline{\Rees}}\right]_{n+1}
= I \left[\omega_{\overline{\Rees}}\right]_n =
I\tau(I^n)\subseteq \tau (I^{n+1}),
\]
using Remark~\ref{remark:tauInOmegan},
\cite[Proposition~3.2]{HyryCoeffIdeals2001},
the induction hypothesis and
Theorem~\ref{theorem:testIdealProps}\eqref{theorem:testIdealProps:prod}.

We now see from Lemma~\ref{lemma:vanishingzerostar} that
\[
\left[ 0^{*}_{\homology^{d+1}_{\overline{\mathfrak{M}}}
(\overline\Rees)} \right]_{-n}=0
\]
for all $n\geq 1$.
Since the $a$-invariant $a(\overline{\Rees})$ is $-1$
(Theorem~\ref{theorem:cmproperties}~\eqref{theorem:cmproperties:ainv}),
it follows that
\[
0^{*}_{\homology^{d+1}_{\overline{\mathfrak{M}}}
(\overline\Rees)}=0.
\]
Hence $\overline\Rees$ is $F$-rational.
\end{proof}

\begin{remarkbox}
%\label{remarkbox:
We now give a minor generalization of
Theorem~\ref{theorem:frationalityifftestideal}. Suppose, more generally,
that $I$ has a reduction $J$ generated by a system of parameters.
Assume the remaining hypothesis on $I$ from 
Theorem~\ref{theorem:frationalityifftestideal}.
Let $ 1 \leq n \leq d-1$.
Then
\[
\tau({I}^n)  = 
\tau({J}^n) \subseteq {J}^n :_{R}\overline{{J}^{d-1}} 
\subseteq {I}^n :_{R}\overline{{I}^{d-1}}
= \tau(I^n).
\]
(The first equality and third inclusion hold since ${J}^n$ is a reduction
of ${I}^n$;
the second is from Remark~\ref{remark:tauInOmegan};
the last equality is by hypothesis.)
Therefore $\tau({J}^n) = {J}^n :_{R}\overline{{J}^{d-1}}$.
Hence $\overline{\Rees(I) } = \overline{\Rees(J) }$ is $F$-rational.
\end{remarkbox}

\section{Corollaries}
\label{section:cor}
\numberwithin{equation}{section}

In this section, we prove some corollaries of the results of the previous
section. Throughout this section, we will assume that $(R, \frakm)$ and $I$
are as in Setup~\ref{setup:main}. We start with the proof of
Proposition~\ref{proposition:gorG}.

\begin{proof}[Proof of Proposition~\protect{\ref{proposition:gorG}}]
\eqref{proposition:gorG:fwd}:
Assume that $\tau(I^{-a-1}) = R$. By Remark~\ref{remark:tauInOmegan},
we need to show that
$\left [\omega_{\overline{\Rees}}\right]_n \subseteq
\tau(I^n)$ for each $n \geq 1$.
Since $\overline{G}$ is Gorenstein and $a<0$, we see that
$\overline{\Rees}$ is Cohen-Macaulay 
(Theorem~\ref{theorem:cmproperties}~\eqref{theorem:cmproperties:cm}),
and that
$\left[ \omega_{\overline{\Rees}}\right]_n = \overline{I^{n+a+1}}$
for all $n\geq 1$ (Proposition~\ref{proposition:descrOmegaR}).
Therefore we will show that
$\overline{I^{n+a+1}}\subseteq \tau (I^n)$ for all $n \geq 1$.
Since
\begin{equation}
\label{equation:tauChain}
\tau (I^{-a-1})\subseteq \tau (I^{-a-2})\subseteq \cdots \subseteq \tau (R)
\end{equation}
we may assume that $n \geq -a$.
In view of~\eqref{equation:tauGor} and
Lemma~\ref{lemma:jlstarnondecreasing}, we will take
$l \gg 0$ and show that
$\overline{I^{n+a+1}} \subseteq I^{[l]} : I^{[l]^{*I^n}}$.

Let $z \in I^{[l]^{*I^n}}$ for some $n> -a-1$ and $l\geq 1$
with $dl-n \geq 1$.
Then there exists a nonzero $c\in R$ such that
$c z^q \in \overline{I^{(dl-n)q}}+I^{[lq]}$ for all $q\gg 1$.
Multiplying by $(\overline{I^{n+a+1}})^{[q]}$, we see that
$c (z\overline{I^{n+a+1}})^{[q]} \subseteq
\overline{I^{(dl+a+1)q}} + I^{[lq]}$ for all $q\gg 1$.
This implies that
$z \overline{I^{n+a+1}} \subseteq
I^{[l]^{*I^{-a-1}}} = \overline{I^{dl+a+1}} + I^{[l]}$, where
the last equality is due to Lemma~\ref{lemma:vanishingzerostar}.
(Note that statement~\eqref{lemma:vanishingzerostar:zerostar} of
Lemma~\ref{lemma:vanishingzerostar} holds with $n=-(a+1)$.)
By Proposition~\ref{proposition:aInvRedNo}
$\overline{I^{dl+a+1}} = I^{dl-d+1} \overline{I^{a+d}}$.
Hence
$z \overline{I^{n+a+1}} \subseteq I^{dl-d+1}  + I^{[l]}
=I^{[l]}$. In other words,
$\overline{I^{n+a+1}} \subseteq I^{[l]}: {I^{[l]^{*I^n}}}$.
This completes the proof of the first assertion of the proposition.

\eqref{proposition:gorG:conv}:
Assume that $\overline{\Rees }$ is $F$-rational and that $a \leq -2$. 
By~\eqref{equation:tauChain} it suffices to show that $\tau(I^{-a-1}) = R$.
Since $\overline{\Rees }$ is $F$-rational, we see
from Remark~\ref{remark:tauInOmegan} and
Proposition~\ref{proposition:descrOmegaR} that
$\tau (I^{-a-1}) = \left[\omega_{\overline{\Rees}} \right]_{-a-1}
= \overline{I^{-a-1+a+1}} = R$.
(Note that, by hypothesis, $\overline{G}$ is Gorenstein and $-a-1 \geq 1$.)
\end{proof}

\begin{corollary}
\label{corollary:frationalinr2}
Suppose that $R$ is a three-dimensional Gorenstein $F$-rational ring
and that $I$ is a reduction of $\frakm$ generated by a system of parameters.
Let $\calI = (\overline{\frakm^n})_{n \in \naturals}$.
Assume that $r_I(\calI) = 2$.
Then $\overline{\Rees}$ is $F$-rational.
\end{corollary}

\begin{proof}
By \cite[Corollary~4.4]{HeinzerKimUlrichCMGor2011} $\overline{G}$ is
Gorenstein.
Further, $a(\overline{G}) = \dim R- r_I(\calI) = -1$.
Now use Proposition~\ref{proposition:gorG}.
\end{proof}

\begin{proof}%
[Proof of Proposition~\protect{\ref{proposition:LevelGaInvMinusOne}}]
Write $\Rees' = R[It, t^{-1}]$ and $\overline{\Rees'}$ for its
normalization. By~\cite[Theorem~1.1]{KoleyKumminiFrationality2017},
We need to show that $\overline{\Rees'}$ is $F$-rational.
Write $\frakM'$ for the homogeneous maximal ideal of $\overline{\Rees'}$.
From the exact sequence
\[
0 \to \overline{\Rees'}(1) \stackrel{t^{-1}}\to \overline{\Rees'}
\to \overline{G} \to 0
\]
we see that for each $n \in \ints$,
\[
\left[\socle
\left( \homology_{\overline{G}_+}^d(\overline{G})\right )\right ]_n
=
\left[\socle
\left( \homology_{\frakM'}^{d+1}(\overline{\Rees'})\right )\right ]_{n+1}.
\]
(Note that $\sqrt{\frakM' \overline{G}} = \sqrt{\overline{G}_+}$ and this
equals the homogeneous maximal ideal of $\overline{G}$.)
It follows that
\[
\socle
\left( \homology_{\frakM'}^{d+1}(\overline{\Rees'})\right )
\subseteq
\left[
\homology_{\frakM'}^{d+1}(\overline{\Rees'})\right ]_{0}.
\]
Arguing as in~\cite[Discussion~3.1]{KoleyKumminiFrationality2017} and using
the fact that $a(\overline{G} ) < 0$, we see that
\[
\left[
\homology_{\frakM'}^{d+1}(\overline{\Rees'})\right ]_{0}
=
\homology_{\frakm}^{d}(R).
\]
We need to show that
\[
0^{*}_{H_{\frakM'}^{d+1}(\overline{\Rees')}}=0.
\]
By way of contradiction, assume that this does not hold. Then there exists
a homogeneous
\[
\xi \in
0^{*}_{H_{\frakM'}^{d+1}(\overline{\Rees')}} \bigcap
\socle
\left( \homology_{\frakM'}^{d+1}(\overline{\Rees'})\right ),
\; \xi \neq 0.
\]
By above, $\deg \xi = 0$. Let $c \in R$ be such that $cF^e(\xi) = 0$ for all
$q \gg 1$.
Then $\xi$ gives a non-zero element of
$0^{*}_{H_{\frakm}^{d}(R)}$, which is a contradiction.
\end{proof}

When $R$ is $F$-rational and $a(\overline{G}) < -1$, then it is not
necessarily true that $\tau(I^{-a-1}) = R$. Indeed, if $R$ and
$\overline{G}$ are Gorenstein and $a(\overline{G}) = -2$, then by
Proposition~\ref{proposition:gorG}, $\tau(I) = R$ if and only if
$\overline{\Rees}$ is $F$-rational.
For a specific example,
see~\cite[Example~6.3]{SinghMultiSymbReesAndStrongFrationality2000}.

We now have the following corollary of
Lemma~\ref{lemma:vanishingzerostar}, relating
the $F$-rationality of Veronese subrings of
$\overline{\Rees}(I)$ to $I^n$-tight closure.
The implication
\eqref{proposition:veronese:proj} $\implies$ \eqref{proposition:veronese:ver}
can also be proved by~\cite[Proposition~5.3]{KoleyKumminiFrationality2017}.

\begin{proposition}
\label{proposition:veronese}
Assume Setup~\ref{setup:main} and that $\overline{\Rees }$ is Cohen-Macaulay.
Then the following are equivalent:
\begin{enumerate}
\item
\label{proposition:veronese:proj}
$\Proj \overline{\Rees}$ is $F$-rational.

\item
\label{proposition:veronese:tc}

For all $n\gg 0$, $I^{[l]^{*I^n}}=\overline{I^{dl-n}}+I^{[l]}$ for all $l$
such that
$dl-n\geq  1$.

\item
\label{proposition:veronese:ver}
For all $n\gg0$, $\overline{\Rees(I^n)}$ is $F$-rational.
\end{enumerate}
\end{proposition}

\begin{proof}

\eqref{proposition:veronese:proj}
$\implies$
\eqref{proposition:veronese:tc}:
Since $\Proj \overline{\Rees}$ is $F$-rational,
$0^{*}_{\homology^{d+1}_{\overline{\mathfrak{M}}} (\overline{\Rees})}$ is
of finite length.
Hence $\left[ 0^{*}_{\homology^{d+1}_{\overline{\mathfrak{M}}}
(\overline{\Rees})} \right]_{-n}=0$ for all $n\gg 0$.
Then by Lemma~\ref{lemma:vanishingzerostar}, for all $n\gg 0$,
$I^{[l]^{*I^n}}=\overline{I^{dl-n}}+I^{[l]}$ for all $l$ such that
$dl-n\geq  1$.

\eqref{proposition:veronese:tc}
$\implies$
\eqref{proposition:veronese:ver}:
Since $\overline{\Rees}$ is Cohen-Macaulay and normal, so is
$\overline{\Rees(I^n)}$ for all $n\geq 1$.
Choose an $n_0$ such that for all $n\geq n_0$,
$I^{[l]^{*I^n}}=\overline{I^{dl-n}}+I^{[l]}$ for all $l$ such that
$dl-n\geq  1$.
By Lemma~\ref{lemma:vanishingzerostar},
$\left [0^{*}_{\homology^{d+1}_{\overline{\mathfrak{M}}}
(\overline{\Rees})} \right]_{-n}=0$ for all $n\geq n_0$.
Since 
\[
\left[
\homology^{d+1}_{\overline{\mathfrak{M}}_{\overline{\Rees(I^{n_0})}}}
(\overline{\Rees(I^{n_0})}) \right]_{-n}=\left[
\homology^{d+1}_{\overline{\mathfrak{M}}_{\overline{\Rees}}}
(\overline{\Rees}) \right]_{-nn_0}
\]
for all $n\geq 1$, it follows that
\[
\left[ 0^{*}_{
\homology^{d+1}_{\overline{\mathfrak{M}}_{\overline{\Rees(I^{n_0})}}}
(\overline{\Rees(I^{n_0})})} \right]_{-n}=0
\]
for all $n\geq 1$. Hence
\[
0^{*}_{
\homology^{d+1}_{\overline{\mathfrak{M}}_{\overline{\Rees(I^{n_0})}}}
(\overline{\Rees(I^{n_0})})}=0. 
\]
Consequently, $\overline{\Rees(I^{n_0})}$
is $F$-rational. Therefore $\overline{\Rees(I^n)}$ is $F$-rational, for all
$n\gg0$.

\eqref{proposition:veronese:ver}
$\implies$
\eqref{proposition:veronese:proj}: Note that
$\Proj\overline{\Rees(I^n)} =\Proj \overline \Rees$.
\end{proof}

\section{Preliminaries, II}
\label{section:prelimsII}
\numberwithin{equation}{subsection}

\subsection{Hypersurface rings}
We now collect some facts about quotients rings of power series rings by a
non-zero power series.

For an element $f$ in a noetherian local ring $(R, \frakm)$, the
\define{order} of $f$, written $\ord(f)$, is
the (unique) integer $k$ such that 
$f \in \frakm^k \minus \frakm^{k+1}$.
We need the following corollary of the Weierstrass Preparation Theorem;
see~\cite[Corollary~9.6]{LeuschkeWiegandCMReprs2012}.

\begin{theorem}
\label{theorem:wptcor} Let $\Bbbk$ be an infinite field and let $f$ be
a non-zero power series in $S=\Bbbk[\![X_0,\cdots X_n]\!]$, $n\geq 1$.
Assume that $\ord(f)=e \geq 2$ and that $e \neq 0 \in \Bbbk$.
Then, after a change of coordinates, we have
$f =u(X_0^e + b_2 X_0^{e-2}+ b_3 X_0^{e-3}+\cdots+ b_{e-1} X_0 + b_e )$,
where $u$ is a unit of $S$ and $b_2,\cdots , b_e$ are non-units of
$\Bbbk[\![X_1,\cdots ,X_n]\!]$.
\end{theorem}

Let $I$ be an $R$-ideal.  An element $r \in I$ is said to be a
\define{superficial element} of $I$ if there exists $c\in
\mathbb{N}$ such that for all $n\geq c$, $(I^{n+1}:r) \cap I^c = I^n$.
A sequence of elements $r_1, r_2,\cdots r_n\in I$ is said to be a
\define{superficial sequence} for $I$ if for all $i=1,2,\cdots, n$ the
image of $x_i$	in $I/(r_1,\cdots, r_{i-1})$ is a superficial element of
$I/(r_1,\cdots,r_{i-1})$.
If $r\in I\minus I^2$
is such that $r^* \in G(I)_1$ is a non-zero-divisor on $G(I)$, then
$r$ is a superficial element of $I$.

\begin{lemma}
\label{lemma:superficialelement}
Let $d \geq 1$.
Let $g \in (X_1,X_2,\cdots,X_d)^2
\subseteq \Bbbk[\![X_1,X_2,\cdots,X_d]\!]$
and $f = X_0^2 + g \in S := \Bbbk[\![X_0,X_1,\cdots,X_d]\!]$.
Write $(R, \frakm ) = S/(f)$.
Then for each $1 \leq i \leq d$, the image of $X_i$ in $R$ is a superficial
element of $\frakm$. In particular, the images of
$X_1, \ldots, X_d$ is a superficial sequence for $\mathfrak{m}$.
\end{lemma}

\begin{proof}
Write $x_0, \ldots, x_d$ for the images of $X_0, \ldots, X_d$ in $R$ and
$Y_0, \ldots, Y_d$ for the images of $x_0, \ldots, x_d \in
\frakm/\frakm^2$.
Write $g(X_1, \ldots, X_d) = g_2(X_1, \ldots, X_d) + g'(X_1, \ldots, X_d)$
with $g'(X_1, \ldots, X_d) \in (X_1,X_2,\cdots,X_d)^3$.
Then
$G(\frakm) \simeq \Bbbk[Y_0, \ldots, Y_d]/(Y_0^2 + g_2(Y_1, \ldots, Y_d))$.
Therefore for all $1 \leq i \leq d$, $Y_i$ is a non-zero-divisor on
$G(\frakm)$, so $x_i$ is a superficial element of $\frakm$.

The second assertion now follows by induction on $d$.
\end{proof}

\subsection{Binomial coefficient modulo a prime}

Let
$m,n$ be two positive integers and $p$  a prime number. Assume that $
m=m_{k}p^{k}+m_{k-1}p^{k-1}+\cdots +m_{1}p+m_{0}$ and
$n=n_{k}p^{k}+n_{k-1}p^{k-1}+\cdots +n_{1}p+n_{0}$ are the base $p$
expansions of $m$ and $n$ respectively, i.e., $m_i,n_i$ are integers such
that $0\leq m_i,n_i \leq p-1$ for all $0\leq i\leq k$.  Then by Lucas's
theorem \cite[Theorem 1]{FineBinomialCoeffs1947},
\begin{equation}
\label{equation:lucas}
\binom {m}{n} \equiv
\prod _{i=0}^{k}{\binom {m_{i}}{n_{i}}} \pmod {p}.
\end{equation}
This uses the convention that $\tbinom {m}{n}=0$ if $m < n$ and that
$\tbinom{0}{0}=1$.
Thus $\binom{m}{n}\not \equiv 0 \pmod{p}$ if and only if $m_i\geq n_i$ for
all $i$.

\begin{lemma} \label{lemma:unitinmodp}
Let $p$ be a prime number.
The following quantities are not divisible by $p$:

\begin{enumerate}

\item
\label{lemma:unitinmodp:pe}
$\binom{p^e-1}{r}$ for each integer  $e\geq 1$ and $0\leq r\leq p^e-1$.

\item
\label{lemma:unitmodp:pe+1}
$\dbinom{\frac{p^e+1}{2}}{1}$ for all integers $e\geq 1$ and $p>2$.

\item
\label{lemma:unitinmodp:beta}
$\dbinom{\frac{p^2+1}{2}}{\beta}$, where $p>3$ and
\[
\beta=
\begin{cases}
\frac{p^2-1}{3} & \text{ if } p\equiv 1 \pmod{3} \\
\frac{2p^2-p+3}{6} & \text{ if } p\equiv 2 \pmod{3}
\end{cases}.
\]

\end{enumerate}
\end{lemma}

\begin{proof}
\eqref{lemma:unitinmodp:pe}:
Since $r\leq p^e-1$, we can write the base-$p$ representation of $r$ as
$r=r_{e-1}p^{e-1}+r_{e-2}p^{e-2}+\cdots +r_0$.
The base-$p$ representation of $p^e-1$is
$p^e-1=(p-1)p^{e-1}+(p-1)p^{e-2}+\cdots +(p-1)$.
Hence, by~\eqref{equation:lucas},
\[
\binom{p^e-1}{r}
\equiv
\prod _{i=0}^{e-1}{\binom {p-1}{r_{i}}}
\not \equiv
0
\pmod {p}.
\]

\eqref{lemma:unitmodp:pe+1}
Note that $\frac{p^e+1}{2}=\frac{p^e-1}{2}+1=\frac{p-1}{2}p^{e-1}+\frac{p-1}{2}p^{e-2}+\cdots +\frac{p-1}{2}+1$. So $\frac{p^e+1}{2}\equiv \frac{p+1}{2}$ $\pmod{p}$. So $\dbinom{\frac{p^e+1}{2}}{1} \not \equiv 0 \pmod{p}$.

\eqref{lemma:unitinmodp:beta}:
Note that the base-$p$ representation of $\frac{p^2+1}{2}$ is
$\frac{p^2+1}{2}=\frac{p-1}{2}p+\frac{p+1}{2}$.
If $p\equiv 1 \pmod{3}$, then the base-$p$ representation of
$\frac{p^2-1}{3}$ is $\frac{p^2-1}{3} = \frac{p-1}{3}p+\frac{p-1}{3}$.
Since $\frac{p+1}{2}> \frac{p-1}{2}> \frac{p-1}{3}$, we see
from~\eqref{equation:lucas} that
\[
\dbinom{\frac{p^2+1}{2}}{\frac{p^2-1}{3}}\not \equiv 0 \pmod{p}.
\]

Now suppose that $p\equiv 2 \pmod{3}$.
The base-$p$ representation of $\frac{2p^2-p+3}{6}$ is
$\frac{2p^2-p+3}{6}=\frac{p-2}{3}p+\frac{p+1}{2}$.

As $\frac{p-1}{2}> \frac{p-2}{3}$, we see  by the above
discussion that
$\dbinom{\frac{p^2+1}{2}}{\frac{2p^2-p+3}{6}}\not \equiv 0 \pmod{p}$.
\end{proof}

\section{Proof of Theorem~\protect{\ref{theorem:projimplies}}}

\label{section:pfThmProjImplies}
\numberwithin{equation}{section}

We now list some sufficient conditions for two-dimensional hypersurface
rings to be $F$-rational.
They are proved using~\cite[Theorem~2.3]{GlassbrStrongFreg1996}.

\begin{lemma}
\label{lemma:hypSurfFRationalTwo}
Let $S = \Bbbk[\![X_0,X_1,X_2]\!]$ and $R=S/(U_0X_0^2+U_1X_1^2+U_2X_2^m)$
where $U_0,U_1,U_2$ are units in $S$ and $m\geq 2$.
Assume that $R$ is $F$-finite and that $p \geq 3$.
Then $R$ is strongly $F$-regular. In particular it is  $F$-rational.
\end{lemma}

\begin{proof}
Write $f = U_0X_0^2+U_1X_1^2+U_2X_2^m$.
Since $X_0$ is in the jacobian ideal of $f$, $R_{X_0}$ is regular.
We want to show that $X_0 f^{p^e-1} \not \in (X_0, X_1, X_2)^{[p^e]}$ for
some $e$, in order to apply~\cite[Theorem~2.3]{GlassbrStrongFreg1996}.
Since $(X_0, X_1, X_2)^{[p^e]}$ is a monomial ideal, it suffices
to exhibit a term of $X_0 f^{p^e-1}$ that does not belong to
$(X_0, X_1, X_2)^{[p^e]}$.

Consider the monomial
$g := X_0(X_0^2)^{\frac{p^{e}-3}{2}}(X_1^2)^{\frac{p^{e}-1}{2}}X_2^m$.
Then for all $e>m$, $g \notin (X_0,X_1,X_2)^{[p^{e}]}$.
Note that the coefficient of $g$ in $X_0 f^{p^e-1}$ is
\[
\dbinom{p^{e}-1}{\frac{p^{e}-3}{2}} \cdot \frac{p^{e}+1}{2}
\]
times a non-zero element of $\Bbbk$.
This is non-zero by Lemma~\ref{lemma:unitinmodp}.
Hence $X_0 f^{p^e-1} \not \in (X_0, X_1, X_2)^{[p^e]}$.
\end{proof}

\begin{lemma}
\label{lemma:hypSurfFRationalThree}
Let $S = \Bbbk[\![X_0,X_1,X_2]\!]$ and
$R=S/(U_0X_0^2+U_1X_1^3+U_2X_1X_2^m +U_3X_2^n)$.
Assume the following:
\begin{enumerate}

\item
$R$ is $F$-finite, $p \geq 7$ and $U_0$ and $U_1$ are invertible elements
of $S$.

\item
For $i = 2,3$, if $U_i$ is non-zero, then it is invertible.
At least one of $\{U_2,U_3\}$ is non-zero.

\item
($U_2\neq 0$ and $2\leq m\leq 3$) or
($U_3\neq 0$ and $3\leq n\leq 5$).

\end{enumerate}

Then $R$ is strongly $F$-regular. In particular it is  $F$-rational.
\end{lemma}

\begin{proof}
Write $f = U_0X_0^2+U_1X_1^3+U_2X_1X_2^m +U_3X_2^n$.
Since $X_0$ is in the jacobian ideal of $f$, $R_{X_0}$ is regular.
We will show that
$X_0 f^{p^2-1} \not \in (X_0, X_1, X_2)^{[p^2]}$
and apply~\cite[Theorem~2.3]{GlassbrStrongFreg1996}.
As in the proof of the previous lemma, it suffices
to exhibit a term of $X_0 f^{p^2-1}$ that does not belong to
$(X_0, X_1, X_2)^{[p^2]}$.
Let
\[
\alpha=\begin{cases}
\frac{p^2-1}{3} & \text{ if } p\equiv 1 \pmod{3}\\
\frac{2p^2-p+3}{6} & \text{ if } p\equiv 2 \pmod{3}
\end{cases}
\]
and $\beta = \frac{p^2+1}{2} - \alpha$.
Note that, by Lemma~\ref{lemma:unitinmodp},
\begin{equation}
\label{equation:nzCoeff}
\dbinom{p^2-1}{\frac{p^2-3}{2}}\dbinom{\frac{p^2+1}{2}}{\alpha}
\in \Bbbk^\times.
\end{equation}

Case I: $U_2 \neq 0$ and $2\leq m\leq 3$.
Consider the monomial
$g := X_0(X_0^2)^{\frac{p^2-3}{2}}(X_1^3)^{\beta}(X_1X_2^m)^{\alpha}
= X_0^{p^2-2}X_1^{3\beta +\alpha} X_2^{m\alpha}$. Then
$g \notin (X_0,X_1,X_2)^{[p^2]}$.
The coefficient of $g$ in $X_0f^{p^2-1}$
is non-zero by~\eqref{equation:nzCoeff}.
Hence $X_0f^{p^2-1 } \not \in (X_0,X_1,X_2)^{[p^2]}$.

Case II: $U_3 \neq 0$ and $3\leq n\leq 5$.
Similar argument as above, with
$g := X_0(X_0^2)^{\frac{p^2-3}{2}}(X_1^3)^{\alpha}(X_2^n)^{\beta} =
X_0^{p^2-2}X_1^{3\alpha}X_2^{n\beta}$.
\end{proof}

We now prove a minor modification
of~\cite[Corollary~5.7]{HaraWatanabeYoshidaFrationality2002}.

\begin{lemma}
\label{lemma:superficialreduction}
Let $(R,\mathfrak{m})$ and $I$ be as in Setup~\ref{setup:main}.
Assume that $\overline{G}$ is Cohen-Macaulay.
Let $x \in \overline{I} \minus \overline{I^2}$ be such that its image in
$\overline{G}$ is a non-zero-divisor on $\overline{G}$.
Write $S = R/(x)$.
If $S$ and the Rees algebra $\overline{S[ISt]}$ are $F$-rational,
then
$R$ and $\overline{\Rees}$ are $F$-rational.
\end{lemma}

\begin{proof}
By~\cite[Theorem~1.1]{KoleyKumminiFrationality2017}, the extended Rees
algebra $\overline{S[ISt,t^{-1}]}$ is
$F$-rational. Since $\overline{G}$ is Cohen-Macaulay,
$\overline{R[It,t^{-1}]}$ is Cohen-Macaulay, and, therefore,
\[
\overline{R[It,t^{-1}]} / (xt) \simeq
\overline{S[ISt,t^{-1}]}.
\]
(It follows from~\cite[Theorem~1]{ItohNormalHilbCoeffs92} that, in general,
$\overline{R[It,t^{-1}]} / (xt)$ and $\overline{S[ISt,t^{-1}]}$ agree
except in finitely many positive degrees.
Since, in our situation, $\depth (\overline{R[It,t^{-1}]} / (xt)) \geq 2$
and $\depth (\overline{S[ISt,t^{-1}]}) \geq 1$, we get the above
isomorphism.)
Hence $\overline{R[It,t^{-1}]}$ is $F$-rational,
by~\cite[Corollary~5.7]{HaraWatanabeYoshidaFrationality2002}).
Hence $R$ and $\overline{\Rees}$ are
$F$-rational~\cite[Theorem~1.1]{KoleyKumminiFrationality2017}.
\end{proof}

\begin{proposition}
\label{proposition:hypdeg2}
Let $S = \Bbbk[\![X_0,\cdots,X_d]\!]$ where $d \geq 2$ and
$\Bbbk$ is an infinite $F$-finite field of characteristic $p \geq 7$.
Let $f \in S$ be irreducible of order $2$; write $R = S/(f)$ and $\frakm$
for the maximal ideal of $R$.
Assume that $\Proj(\Rees(\frakm))$ is $F$-rational.
Then $R$ and $\Rees(\frakm)$ are $F$-rational.
\end{proposition}

\begin{proof}
The strategy of the proof is to show (using the $F$-rationality of $\Proj
\Rees(\frakm)$) that after a suitable change of
variables, $f$ is of the form $\hat{f} + g$ where $\hat{f}$ the power
series in Lemmas~\ref{lemma:hypSurfFRationalTwo}
or~\ref{lemma:hypSurfFRationalThree} and $g \in (X_3, \ldots, X_d)S$.
Hence $R/(X_3, \ldots, X_d)R \simeq \Bbbk[\![X_0, X_1, X_2]\!]/(\hat{f})$
is a two-dimensional $F$-rational ring.
By~\cite[Theorem 3.1]{HaraWatanabeYoshidaFrationality2002}, the Rees
algebra of an integrally closed ideal in a
two-dimensional $F$-rational ring is $F$-rational.
Now apply Lemma~\ref{lemma:superficialreduction} repeatedly.

We first observe that $\Rees(\frakm)$ is Cohen-Macaulay and normal.
Since $G(\frakm )$ is Gorenstein with $a$-invariant $2-(d+1)$,
$\Rees(\frakm)$ is Cohen-Macaulay.
By Remark~\ref{remarkbox:reesnormal}, it is normal.

By Theorem~\ref{theorem:wptcor}, we may assume that $f=X_0^2+g$,
where $g \in \Bbbk[\![X_1, \ldots, X_d]\!]$ and $\ord(g) \geq 2$.
(We are concerned only about the ideal $fS$, not the element $f$,
\textit{per se}.)
Since $R$ is a domain, $g\neq 0$.

Let $I=(X_1,\cdots,X_d)R$.
Then $I$ is a minimal reduction of $\frakm$ and $I\frakm = \frakm^2$.
By Lemma~\ref{lemma:superficialelement}, $x_1,\cdots, x_d$ is a sequence of
superficial elements.

We now show that $\ord(g)\leq 3$.
By way of contradiction, assume that $\ord(g)\geq 4$.
Then for all $q \geq 1$,
$X_0^{2q} \equiv g^q \mod f$, so in $R$, $X_0^{2q} \in I^{4q}$.
Hence $X_0 \in \overline{I^2} = \overline{\frakm^2} = \frakm^2$,
by~\cite[6.8.3]{SwHuIntCl06} and the normality of $\Rees(\frakm)$.
This is a contradiction. Hence $\ord(g)\leq 3$.

Case I: $\ord(g)=2$.
By Theorem~\ref{theorem:wptcor}, we may assume that
$g=U_1(X_1^2+h)$, where $U_1$ is a unit in
$\Bbbk[\![X_1,\cdots,X_d]\!]$ and
$h \in \Bbbk[\![X_2,\cdots, X_d]\!]$ with $\ord(h)\geq 2$.

We now show that $h \neq 0$.
By way of contradiction, assume that $h=0$.
Then $f = X_0^2 + U_1X_1^2 \in (X_0, X_1)S$.
The ring
$R_{(X_0, X_1)}$ is not normal, since it is one-dimensional, but not
regular.
($f$ is in the square of the maximal ideal of $S_{(X_0, X_1)}$.)
Hence $\Spec R \minus \{ \frakm\}$ is not normal.
In particular, $\Proj \Rees(\frakm)$ is not normal, which contradicts the
hypothesis that $\Proj \Rees(\frakm)$ is $F$-rational.
Hence $h \neq 0$.

Let $m = \ord(h)$. Again by Theorem~\ref{theorem:wptcor},
write $h=U_2(X_2^m+h')$, where $U_2$ is a unit in
$\Bbbk[\![X_2,\cdots,X_d]\!]$ and
$h'\in (X_3,\cdots,X_d)\Bbbk[\![X_2,\cdots,X_d]\!]$
with $\ord (h')\geq \ord (h)$.
Let
\[
\overline{R} := \frac{R}{(X_3,\cdots,X_d)} \simeq
\frac{\Bbbk[\![X_0,X_1,X_2]\!]}{(X_0^2+U_1X_1^2+U_2X_2^m)}
\]
for some units $U_1,U_2\in \Bbbk[\![X_0,X_1,X_2]\!]$.
By Lemma~\ref{lemma:hypSurfFRationalTwo}, $\overline{R}$ is $F$-rational.
By~\cite[Theorem~3.1]{HaraWatanabeYoshidaFrationality2002}, the Rees
algebra $\overline{R}[X_0t, X_1t, X_2t]$ is $F$-rational. By
Lemmas~\ref{lemma:superficialelement} and~\ref{lemma:superficialreduction},
$\Rees(\frakm)$ is $F$-rational.

Case II: $\ord(g)=3$.
By Theorem~\ref{theorem:wptcor}, we may assume that
$g=U_1(X_1^3+X_1h_1+h_2)$, where $U_1$ is a unit in $\Bbbk[\![X_1,\cdots,
X_d]\!]$ and $h_1,h_2$ are non-units in $\Bbbk[\![X_2,\cdots,X_d]\!]$ such
that $m := \ord(h_1)\geq 2$ and $n := \ord(h_2)\geq 3$.
As in the above case, if $h_1 = h_2 = 0$, then $f$ is in the square of the
maximal ideal of $S_{(X_0,X_1)}$, so
$\Spec R \minus \{ \frakm\}$ would not be normal, contradicting the
$F$-rationality of $\Proj \Rees(\frakm)$. Hence $h_1 \neq 0$ or
$h_2 \neq 0$.

Note that $(X_0)^{q+1}X_1^{lq}=(-g)^{\frac{q+1}{2}}X_1^{lq}=
\pm (U_1X_1^3+h_1X_1+h_2)^{\frac{q+1}{2}}X_1^{lq}$
and that a general term of the above expression is of the form
$X_1^{3\alpha+\beta +lq}h_1^{\beta}h_2^{\gamma}$, with some coefficient,
where $\alpha+\beta+\gamma=\frac{q+1}{2}$.

Now suppose that $m\geq 4$ and $n\geq 6$. We claim that
\begin{equation}
\label{equation:term}
X_1^{3\alpha+\beta +lq}h_1^{\beta}h_2^{\gamma} \in I^{(l+2)q}+I^{[(l+1)q]}
\end{equation}
for all $q \gg 1$ and $l \gg 0$.
Assume the claim.
Then
\[
X_0 (X_0X_1^l)^q \in I^{(l+2)q}+I^{[(l+1)q]}
\]
for all $q \gg 1$ and $l \gg 0$.
By~\eqref{equation:jlstardefn}
$X_0X_1^l\in I^{[l+1]*^{I^{(d-1)(l+1)-1}}}$ for all $l \gg 0$.
Since $\Proj \Rees (\frakm )$ is $F$-rational,
\[
\left[ 0^{*}_{\homology^{d+1}_{\mathfrak{M}} (\Rees)}
\right]_{-((d-1)(l+1)-1)} = 0
\]
for all $l \gg 0$, so by Lemma~\ref{lemma:vanishingzerostar}, we see that
$X_0X_1^l\in \frakm^{l+2} + I^{[(l+1)]}$
for all $l \gg 0$.
Then $x_0(x_1x_2\cdots x_d)^l\in \mathfrak{m}^{dl+2}+I^{[l+1]}$. So $x_0\in
(\mathfrak{m}^{dl+2}+I^{[l+1]}):(x_1x_2\cdots x_d)^l\subseteq
\mathfrak{m}^2+I\subseteq I$, which is a contradiction. (Use an argument
as in the proof of Proposition~\ref{proposition:imagePhi}.)
This is a contradiction.
Hence $m\leq 3$ or $n \leq 5$. Using Theorem~\ref{theorem:wptcor}, we may
write $g = U_1(X_1^3+U_2X_1X_2^m+U_3X_2^n+H_3)$ for some
$H_3\in (X_3,\cdots,X_d)\Bbbk[\![X_1,\cdots,X_d]\!]$ and $U_2,U_3$ that are
either invertible or zero but not zero simultaneously.
Let
\[
\overline{R} := \frac{R}{(X_3,\cdots, X_d)} \simeq
\frac{\Bbbk[\![X_0,X_1,X_2]\!]}
{(X_0^2+U_1X_1^3+U_1U_2X_1X_2^{m}+U_1U_3X_2^n)}.
\]
By Lemma~\ref{lemma:hypSurfFRationalThree}, $\overline{R}$ is $F$-rational.
As in the earlier case,
use~\cite[Theorem~3.1]{HaraWatanabeYoshidaFrationality2002} and
Lemmas~\ref{lemma:superficialelement} and~\ref{lemma:superficialreduction}.

It remains to prove the claim~\eqref{equation:term}.
If $ 3\alpha+\beta \leq q-1$ and $3\alpha+\beta +m \beta +n \gamma \leq
2q-1$, then $(q+4-(3\alpha+\beta))+(m-4)\beta +(n-6)\gamma \leq 0$; which
gives a contradiction. Hence either $ 3\alpha+\beta \geq q$ or
$3\alpha+\beta +m \beta +n \gamma \geq 2q$. In either case,
$X_1^{3\alpha+\beta +lq}h_1^{\beta}h_2^{\gamma} \in I^{(l+2)q}+I^{[(l+1)q]}$.
\end{proof}

\begin{proposition}
\label{proposition:dim3}
Let $(R,\mathfrak{m})$ be a three-dimensional Gorenstein $F$-finite
$F$-rational complete local domain of characteristic $p\geq 7$ with an
infinite residue field.
Suppose that $\Proj \Rees(\mathfrak{m})$ is $F$-rational and
that $\Rees(\mathfrak{m})$ is Cohen-Macaulay. Then
$\Rees(\mathfrak{m})$ is $F$-rational.
\end{proposition}

\begin{proof}
Since $\Rees(\mathfrak{m})$ is Cohen-Macaulay, the reduction number
$r(\mathfrak{m})$ is at most $2$.

Case I: $r(\mathfrak{m})=0$. Then $R$ is a regular local ring; hence
$\Rees(\mathfrak{m})$ is $F$-rational.

Case II: $r(\mathfrak{m})=1$.
Then $\Rees(\mathfrak{m})$ is
Gorenstein~\cite[Theorem 4.4]{HermannRibbeZarzuelaGorRees1994}.
Consequently $R$ is a hypersurface with multiplicity
$e(R)=2$~\cite[Corollary 4.5]{HermannRibbeZarzuelaGorRees1994}.
Now use Proposition~\ref{proposition:hypdeg2}.

Case III: $r(\mathfrak{m})=2$.
As we observed in Remark~\ref{remarkbox:reesnormal}
$\Rees(\mathfrak{m})$ is normal. Now apply
Corollary~\ref{corollary:frationalinr2}.
\end{proof}

\begin{proof}[Proof of Theorem~\protect{\ref{theorem:projimplies}}]

\eqref{theorem:projimplies:hyp}:
Follows from Proposition~\ref{proposition:hypdeg2}.

\eqref{theorem:projimplies:dim3}:
Follows from Proposition~\ref{proposition:dim3}.
\end{proof}

It is not difficult to find examples of non-$F$-pure (even
non-$F$-injective) hypersurfaces
$(R, \frakm)$ such that $\Proj \Rees(\frakm)$ is $F$-rational but
$\Rees(\frakm)$ is not.
E.g., look at generic homogeneous polynomials in $d$ variables with degree
at least $d$.
Below, we give an example of an $F$-pure hypersurface $(R, \frakm)$ such
that $\Proj \Rees(\frakm)$ is $F$-rational but $\Rees(\frakm)$ is
Cohen-Macaulay, normal but not $F$-rational or $F$-injective.
It also shows that some hypothesis on the characteristic is necessary in
Proposition~\ref{proposition:hypdeg2}.

\begin{example}
Let $\Bbbk$ is an algebraically closed field of characteristic $2$,
$S = \Bbbk[\![X,Y,Z,W]\!]$, $f = X^2+XYZW+Y^3+Z^3+W^3$ and
$R = \frac{S}{(f)}$.
Write $\mathfrak{m}=(X,Y,Z,W)R$ and $I=(Y,Z,W)R$.
Then we have the following:

\begin{enumerate}

\item
The singular locus of $R$ is $\{\frakm \}$.
Since $f \not \in (X^2, Y^2, Z^2, W^2)S$, $R$ is $F$-pure
by~\cite[Proposition~1.7]{FedderFpurity1983}.
On the other hand, $X^{q}\in (XYZW)^{\frac{q}{2}}R+ I^{[q]}$ for all $q\geq
2$.
Apply induction on $q$ to show that
\[
XX^q \in XX^{\frac{q}{2}} (YZW)^{\frac{q}{2}}R+I^{[q]} \subseteq
(YZW)^{[\frac{q}{2}]} I^{[\frac{q}{2}]}+I^{[q]}\subseteq I^{[q]}
\]
for all $q\geq 2$.
(The second inclusion is by inductive hypothesis.)
Hence $X\in I^{*} \minus I$, so $R$ is not $F$-rational.

\item The associated graded ring $G := G(\mathfrak{m})$ is
Gorenstein; its $a$-invariant is $-2$.
Hence $\Rees := \Rees(\frakm)$ is Gorenstein, and, in particular, Cohen
Macaulay.
The reduction number $r_I(\frakm)$ is $1$, i.e.,
$I\mathfrak{m}^n=\mathfrak{m}^{n+1}$ for each $n \geq 1$.

\item We now show that $\Proj\Rees$ is $F$-rational.
For this, it suffices to show that $\Rees$ is $F$-rational
except at its homogeneous maximal ideal $\frakM$.
Since $\Rees$ is Cohen-Macaulay, it would then follow that $\Rees$ is
normal.

\begin{enumerate}

\item
For all $a \in \frakm$, $\Rees_a \simeq R_a[t]$ is regular.

\item
$\mathscr{R}_{yt} \cong
\frac{\Bbbk[X,Y,Z,W]}{(X^2+Y(XYZW+1+Z^3+W^3))}$. Its non-regular
locus is defined by $(X,Y,1+Z^3+W^3)$.
Choose $\alpha,\beta \in \Bbbk$ be such that $1+\alpha^3+\beta^3=0$.
It is enough to show localization of $\mathscr{R}_{yt}$ at
$(X,Y,Z-\alpha,W-\beta)$ is $F$-rational.

Replacing $Z$ by $Z+\alpha$ and $W$ by $W+\beta$, it is enough to show
that $A := \frac{\Bbbk[X,Y,Z,W]}{g}$ is $F$-rational, where
\begin{align*}
g & = X^2+Y(XY(Z+\alpha) (W+\beta)+1+(Z+\alpha)^3+(W+\beta)^3)
\\
& = X^2+Y(\alpha^2Z+\beta^2 W)+XY^2(Z+\alpha)(W+\beta) + YZ^2(Z+\alpha)
+ YW^2(W+\beta).
\end{align*}
Let $I=(Y,Z,W)A$; it is a minimal reduction of $(X,Y,Z,W)A$, with
$I(X,Y,Z,W)^nA = (X,Y,Z,W)^{n+1}A$ for each $n \geq 1$.
Hence it is enough to show $X\not \in I^{*}$.
By way of contradiction, assume that $X\in I^{*}$.
Then there exist a nonzero element $c$ such that $cX^q\in I^{[q]}$
for all $q\gg 1$.
Note that $cX^q = c(X^2)^{\frac{q}{2}} \equiv cY^{\frac{q}{2}}(\alpha^q
Z^{q/2}+\beta^q W^{q/2}) \mod I^{[q]}$ for all $q\gg 1$.
Hence $cY^{\frac{q}{2}}(\alpha^q Z^{q/2}+\beta^q W^{q/2})\in I^{[q]}$.
Since $Y,Z,W$ is a regular sequence,  we see that
$c(\alpha^q Z^{q/2}+\beta^q W^{q/2}) \in (Y^{\frac{q}{2}},Z^q,W^q)$
for all $q\gg 1$.
Since $\alpha^3 + \beta^3 = 1$, assume, without loss of generality,
that $\alpha\neq 0$.
Use a similar argument to see that $c\alpha^q \in I^{[\frac{q}{2}]}$
for all $q\gg1 $, which is a contradiction.
Hence $\Rees_{yt}$ is $F$-rational.

\item
Since $\Rees_{yt}\cong \Rees_{zt}\cong \Rees_{wt}$, all
these rings are $F$-rational.
Since $\frakM = \sqrt { \frakm \Rees + It\Rees}$ it follows that $\Spec
\Rees \minus \frakM$ is $F$-rational.

\end{enumerate}

\item
Since $R$ is not $F$-rational, $\Rees$ is not $F$-rational
by~\cite[Corollary~2.13]{HaraWatanabeYoshidaFrationality2002}.
One can also see it using Remark~\ref{remark:tauInOmegan} and the
inclusions
\[
\tau(\frakm ) \subseteq \tau(R) \subsetneq R =
[\omega_{\Rees}]_1.
\]
(The first one is a standard property of test ideals, the second holds
since $R$ is not $F$-regular, the final one holds since $\Rees$ is
Cohen-Macaulay and $a(G) = -2$.)

\item We now observe that $\Rees$ is not $F$-injective, and,
therefore, not $F$-pure.
By way of contradiction, suppose that $\Rees$ is $F$-injective.
Since $\Proj \Rees$ is $F$-rational and $\Rees$ is not $F$-rational,
$M := 0^{*}_{\homology^{d+1}_{\frakM}(\Rees)}$
is a non-zero module of finite length.
Let $\xi \in M$
be a non-zero element of minimum degree.
Since $\Rees$ is Cohen-Macaulay, $\deg \xi < 0$.
By the $F$-injectivity of $\Rees$, we see that $\xi^p \neq 0$.
Since $M$ is closed under the application of the Frobenius map, we get a
contradiction of the minimality of $\deg \xi$.
Hence $\Rees$ is not $F$-injective.

\end{enumerate}

\end{example}

\ifreadkumminibib
\bibliographystyle{alphabbr}
\bibliography{kummini}

\begin{thebibliography}{HWY02}

\bibitem[Fed83]{FedderFpurity1983}
R.~Fedder.
\newblock {$F$}-purity and rational singularity.
\newblock {\em Trans. Amer. Math. Soc.}, 278(2):461--480, 1983.

\bibitem[Fin47]{FineBinomialCoeffs1947}
N.~J. Fine.
\newblock Binomial coefficients modulo a prime.
\newblock {\em Amer. Math. Monthly}, 54:589--592, 1947.

\bibitem[Gla96]{GlassbrStrongFreg1996}
D.~Glassbrenner.
\newblock Strong {$F$}-regularity in images of regular rings.
\newblock {\em Proc. Amer. Math. Soc.}, 124(2):345--353, 1996.

\bibitem[GN94]{GotoNishidaFiltrGorRees1994}
S.~Goto and K.~Nishida.
\newblock Filtrations and the {G}orenstein property of the associated {R}ees
  algebras.
\newblock In {\em The Cohen-Macaulay and Gorenstein Rees algebras associated to
  filtrations}, volume 110, pages 69--134. American Mathematical Society,
  Providence, RI, 1994.

\bibitem[HH90]{HochsterHunekeTCInvThyBSThm90}
M.~Hochster and C.~Huneke.
\newblock Tight closure, invariant theory, and the {B}rian\c{c}on-{S}koda
  theorem.
\newblock {\em J. Amer. Math. Soc.}, 3(1):31--116, 1990.

\bibitem[HH94]{HochsterHunekeFregSmBsch94}
M.~Hochster and C.~Huneke.
\newblock {$F$}-regularity, test elements, and smooth base change.
\newblock {\em Trans. Amer. Math. Soc.}, 346(1):1--62, 1994.

\bibitem[HKU11]{HeinzerKimUlrichCMGor2011}
W.~Heinzer, M.-K. Kim, and B.~Ulrich.
\newblock The {C}ohen-{M}acaulay and {G}orenstein properties of rings
  associated to filtrations.
\newblock {\em Comm. Algebra}, 39(10):3547--3580, 2011.

\bibitem[HRZ94]{HermannRibbeZarzuelaGorRees1994}
M.~Herrmann, J.~Ribbe, and S.~Zarzuela.
\newblock On the {G}orenstein property of {R}ees and form rings of powers of
  ideals.
\newblock {\em Trans. Amer. Math. Soc.}, 342(2):631--643, 1994.

\bibitem[HS06]{SwHuIntCl06}
C.~Huneke and I.~Swanson.
\newblock {\em Integral closure of ideals, rings, and modules}, volume 336 of
  {\em London Mathematical Society Lecture Note Series}.
\newblock Cambridge University Press, Cambridge, 2006.

\bibitem[HWY02]{HaraWatanabeYoshidaFrationality2002}
N.~Hara, K.-i. Watanabe, and K.-i. Yoshida.
\newblock F-rationality of {R}ees algebras.
\newblock {\em J. Algebra}, 247(1):153--190, 2002.

\bibitem[HY03]{HaraYoshidaGenTightClMultIdeals2003}
N.~Hara and K.-I. Yoshida.
\newblock A generalization of tight closure and multiplier ideals.
\newblock {\em Trans. Amer. Math. Soc.}, 355(8):3143--3174, 2003.

\bibitem[Hyr99]{HyryBlowupRingsRationalSings1999}
E.~Hyry.
\newblock Blow-up rings and rational singularities.
\newblock {\em Manuscripta Math.}, 98(3):377--390, 1999.

\bibitem[Hyr01]{HyryCoeffIdeals2001}
E.~Hyry.
\newblock Coefficient ideals and the {C}ohen-{M}acaulay property of {R}ees
  algebras.
\newblock {\em Proc. Amer. Math. Soc.}, 129(5):1299--1308, 2001.

\bibitem[Ito92]{ItohNormalHilbCoeffs92}
S.~Itoh.
\newblock Coefficients of normal {H}ilbert polynomials.
\newblock {\em J. Algebra}, 150(1):101--117, 1992.

\bibitem[KK21]{KoleyKumminiFrationality2017}
M.~Koley and M.~Kummini.
\newblock {$F$}-rationality of {R}ees algebras.
\newblock {\em J. Algebra}, 571:151--167, 2021.
\newblock arXiv:1803.05382 [math.AC].

\bibitem[LT81]{LipTeiPseudoRatlSing81}
J.~Lipman and B.~Teissier.
\newblock Pseudorational local rings and a theorem of {B}rian\c con-{S}koda
  about integral closures of ideals.
\newblock {\em Michigan Math. J.}, 28(1):97--116, 1981.

\bibitem[LW12]{LeuschkeWiegandCMReprs2012}
G.~J. Leuschke and R.~Wiegand.
\newblock {\em Cohen-{M}acaulay representations}, volume 181 of {\em
  Mathematical Surveys and Monographs}.
\newblock American Mathematical Society, Providence, RI, 2012.


\bibitem[Sin00]{SinghMultiSymbReesAndStrongFrationality2000}
A.~K. Singh.
\newblock Multi-symbolic {R}ees algebras and strong {F}-regularity.
\newblock {\em Math. Z.}, 235(2):335--344, 2000.

\bibitem[Smi97]{SmithFrational1997}
K.~E. Smith.
\newblock {$F$}-rational rings have rational singularities.
\newblock {\em Amer. J. Math.}, 119(1):159--180, 1997.

\bibitem[Tru87]{TrungRedExpDegBd87}
N.~V. Trung.
\newblock Reduction exponent and degree bound for the defining equations of
  graded rings.
\newblock {\em Proc. Amer. Math. Soc.}, 101(2):229--236, 1987.

\bibitem[V{\'e}l95]{VelezOpennessFrationalLoci1995}
J.~D. V{\'e}lez.
\newblock Openness of the {F}-rational locus and smooth base change.
\newblock {\em J. Algebra}, 172(2):425--453, 1995.

\bibitem[VV78]{ValabregaVallaFormRings1978}
P.~Valabrega and G.~Valla.
\newblock Form rings and regular sequences.
\newblock {\em Nagoya Math. J.}, 72:93--101, 1978.

\end{thebibliography}
\else

\def\cfudot#1{\ifmmode\setbox7\hbox{$\accent"5E#1$}\else
  \setbox7\hbox{\accent"5E#1}\penalty 10000\relax\fi\raise 1\ht7
  \hbox{\raise.1ex\hbox to 1\wd7{\hss.\hss}}\penalty 10000 \hskip-1\wd7\penalty
  10000\box7}

\fi %
\end{document}